\documentclass[final,leqno,onefignum,onetabnum]{siamltex1213}
\usepackage{amsfonts}
\usepackage{amsmath,booktabs,ctable,threeparttable}
\usepackage{amssymb,amsfonts,boxedminipage}
\newtheorem{remark}{Remark}[section]

\newtheorem{algorithm}{Algorithm}[section]

\title{The Krylov Subspaces, Low Rank Approximations and Ritz Values
of LSQR for Linear Discrete Ill-Posed Problems: the Multiple Singular
Value Case\thanks{This work was supported in part by
the National Science Foundation of China (No. 11771249)}}

\author{Zhongxiao Jia\thanks{Department of Mathematical Sciences, Tsinghua
University, 100084 Beijing, China. (\email{jiazx@tsinghua.edu.cn})}}

\begin{document}
\maketitle
\slugger{sirev}{xxxx}{xx}{x}{x--x}

\begin{abstract}
For the large-scale linear discrete ill-posed problem $\min\|Ax-b\|$ or $Ax=b$
with $b$ contaminated by white noise, the Golub-Kahan bidiagonalization based
LSQR method and its mathematically equivalent CGLS, the Conjugate Gradient (CG) method
applied to $A^TAx=A^Tb$, are most commonly used. They have intrinsic regularizing
effects, where the iteration number $k$ plays the role
of regularization parameter. The long-standing fundamental question is:
{\em Can LSQR and CGLS find 2-norm filtering best possible regularized solutions}?
The author has given definitive answers to this question for severely and
moderately ill-posed problems when the singular values of $A$ are simple.
This paper extends the results to the multiple singular
value case, and studies the approximation accuracy of Krylov subspaces, the quality
of low rank approximations generated by Golub-Kahan bidiagonalization and the convergence
properties of Ritz values. For the two kinds of problems, we prove that LSQR finds
2-norm filtering best possible regularized solutions at semi-convergence.
Particularly, we consider some important and untouched issues on best,
near best and general rank $k$ approximations to $A$ for the
ill-posed problems with the singular values $\sigma_k=\mathcal{O}(k^{-\alpha})$
with $\alpha>0$, and the relationships between them and their nonzero
singular values. Numerical experiments confirm our theory.
The results on general rank $k$ approximations and the properties of
their nonzero singular values apply to several
Krylov solvers, including LSQR, CGME, MINRES, MR-II, GMRES and RRGMRES.
\end{abstract}


\pagestyle{myheadings}
\thispagestyle{plain}
\markboth{ZHONGXIAO JIA}{REGULARIZATION OF LSQR IN THE MULTIPLE SINGULAR VALUE CASE}

\section{Introduction and Preliminaries}\label{intro}

Consider the linear discrete ill-posed problem
\begin{equation}
  \min\limits_{x\in \mathbb{R}^{n}}\|Ax-b\| \mbox{\,\ or \ $Ax=b$,}
  \ \ \ A\in \mathbb{R}^{m\times n}, \label{eq1}
  \ b\in \mathbb{R}^{m},
\end{equation}
where the norm $\|\cdot\|$ is the 2-norm of a vector or matrix, and
$A$ is extremely ill conditioned with its singular values decaying
to zero without a noticeable gap. Without loss of generality, we assume
that $m\geq n$. Problem \eqref{eq1} typically arises
from the discretization of the first kind Fredholm integral equation
\begin{equation}\label{eq2}
Kx=(Kx)(s)=\int_{\Omega} k(s,t)x(t)dt=g(s)=g,\ s\in \Omega
\subset\mathbb{R}^q,
\end{equation}
where the kernel $k(s,t)\in L^2({\Omega\times\Omega})$ and
$g(s)$ are known functions, while $x(t)$ is the
unknown function to be sought. If $k(s,t)$ is non-degenerate
and $g(s)$ satisfies the Picard condition, there exists the unique square
integrable solution
$x(t)$; see \cite{engl00,hansen98,hansen10,kirsch,mueller}. Here for brevity
we assume that $s$ and $t$ belong to the same set $\Omega\subset
\mathbb{R}^q$ with $q\geq 1$.
Applications include image deblurring, signal processing, geophysics,
computerized tomography, heat propagation, biomedical and optical imaging,
groundwater modeling, and many others; see, e.g.,
\cite{aster,engl93,engl00,hansen10,kaipio,kirsch,mueller,natterer,vogel02}.
The right-hand side $b=b_{true}+e$ is assumed to be
contaminated by a Gaussian white noise $e$, caused by measurement, modeling
or discretization errors, where $b_{true}$
is noise-free and $\|e\|<\|b_{true}\|$.
Because of the presence of noise $e$ and the extreme
ill-conditioning of $A$, the naive
solution $x_{naive}=A^{\dagger}b$ of \eqref{eq1} bears no relation to
the true solution $x_{true}=A^{\dagger}b_{true}$, where
$\dagger$ denotes the Moore-Penrose inverse of a matrix.
Therefore, one has to use regularization to extract a
best possible approximation to $x_{true}$.

For a Gaussian white noise $e$, we always assume that $b_{true}$ satisfies the
discrete
Picard condition $\|A^{\dagger}b_{true}\|\leq C$ with some constant $C$ for $n$
arbitrarily large \cite{aster,hansen90b,hansen98,hansen10,kern}.
It is an analog of the Picard condition in the finite dimensional case;
see, e.g., \cite[p.9]{hansen98},
\cite[p.12]{hansen10} and \cite[p.63]{kern}.
Without loss of generality, assume that $Ax_{true}=b_{true}$.
Then the two dominating regularization approaches are
to solve the following two equivalent problems:
\begin{equation}\label{posed}
\min\limits_{x\in \mathbb{R}^{n}}\|Lx\| \ \ \mbox{subject to}\ \
\|Ax-b\|\leq \tau\|e\|
\end{equation}
with $\tau\approx 1$ and general-form Tikhonov regularization
\begin{equation}\label{tikhonov}
  \min\limits_{x\in \mathbb{R}^{n}}\{\|Ax-b\|^2+\lambda^2\|Lx\|^2\}
\end{equation}
with $\lambda>0$ the regularization parameter \cite{hansen98,hansen10},
where $L$ is a regularization matrix and its suitable choice is based on
a-prior information on $x_{true}$. Typically, $L$ is either the identity matrix
$I$ or the scaled discrete approximation of a first or second order derivative
operator. If $L=I$, \eqref{tikhonov} is standard-form Tikhonov regularization,
and both \eqref{posed} and \eqref{tikhonov} are 2-norm filtering
regularization problems.

We are concerned with the case $L=I$ in this paper. If $L\not=I$,
\eqref{posed} and \eqref{tikhonov}, in principle,
can be transformed into standard-form
problems \cite{hansen98,hansen10}. In this case, for \eqref{posed} of small or
moderate size, an effective and reliable solution method
is the truncated singular value decomposition (TSVD)
method, and it obtains the 2-norm filtering best regularized solution $x_{k_0}^{tsvd}$
at some $k_0<n$  \cite{hansen98,hansen10}, where
$k_0$ is the optimal regularization parameter, called the transition point, such
that $\|x_{k_0}^{tsvd}-x_{true}\|=\min_{k=1,2,\ldots,n}\|x_k^{tsvd}-x_{true}\|$.
We will review the TSVD method and
reformulate it and \eqref{tikhonov} when $A$ has multiple singular values.
A key of solving \eqref{tikhonov}
is the determination of the optimal regularization parameter $\lambda_{opt}$
such that $\|x_{\lambda_{opt}}-x_{true}\|=\min_{\lambda>0}\|x_{\lambda}-x_{true}\|$.
A number of parameter-choice methods have been developed for finding
$\lambda_{opt}$, such as the discrepancy principle,
the L-curve criterion, and the generalized cross validation
(GCV), etc. We refer the reader to, e.g., \cite{hansen98,hansen10} for details.

It has been theoretically and numerically
justified that $x_{k_0}^{tsvd}$ and $x_{\lambda_{opt}}$
essentially have the minimum 2-norm error;
see \cite{varah79}, \cite{hansen90b}, \cite[p.109-11]{hansen98} and
\cite[Sections 4.2 and 4.4]{hansen10}.
In effect, the theory in \cite{engl00} has shown that the error of
$x_{k_0}^{tsvd}$ is unconditionally order optimal in the Hilbert
space setting, i.e.,
the same order as the worst case error, while $x_{\lambda_{opt}}$
is conditionally order optimal. As a result, we can naturally take
$x_{k_0}^{tsvd}$ as a reference standard when assessing the regularization
ability of a 2-norm filtering regularization method.

For \eqref{eq1} large, the TSVD method and the Tikhonov regularization
method are generally too demanding, and only iterative regularization
methods are computationally viable.
Krylov iterative solvers are a major class of methods for solving \eqref{eq1},
and they project problem \eqref{eq1} onto a sequence of
low dimensional Krylov subspaces
and computes iterates to approximate $x_{true}$
\cite{aster,engl00,gilyazov,hanke95,hansen98,hansen10,kirsch}.
Of them, the CGLS method, which implicitly applies the CG
method \cite{hestenes} to $A^TAx=A^Tb$,
and its mathematically equivalent LSQR algorithm \cite{paige82}
have been most commonly used. The Krylov solvers CGME
\cite{bjorck96,bjorck15,craig,hanke95,hanke01} and
LSMR \cite{bjorck15,fong} are also choices. These Krylov solvers
have general regularizing
effects \cite{aster,gilyazov,hanke95,hanke01,hansen98,hansen10,hps16,hps09}
and exhibit semi-convergence \cite[p.89]{natterer};
see also \cite[p.314]{bjorck96},
\cite[p.135]{hansen98} and \cite[p.110]{hansen10}: The iterates
converge to $x_{true}$ in an initial stage; afterwards the
noise $e$ starts to deteriorate the iterates so that they start to diverge
from $x_{true}$ and instead converge to $x_{naive}$.
If we stop at the right time, then, in principle,
we have a regularization method, where the iteration number plays the
role of the regularization parameter.
Semi-convergence is due to the fact that the projected problem starts to
inherit the ill-conditioning of \eqref{eq1} from some iteration
onwards, and the appearance of a small singular
value of the projected problem amplifies the noise considerably.

The behavior of \eqref{eq1}, \eqref{posed} and \eqref{tikhonov}
with $L=I$ critically depends on the decay rate of
the singular values $\sigma_j$ of $A$.
The behavior of ill-posed problems critically depends on the decay rate of
$\sigma_j$. For a linear compact operator equation
such as (\ref{eq2}) in the Hilbert space setting, let $\mu_1\geq\mu_2\geq\cdots\geq 0$
be the singular values of the compact operator $K$.
The following characterization of the degree of ill-posedness
of (\ref{eq2}) was introduced in \cite{hofmann86}
and has been widely used; see, e.g., \cite{aster,engl00,hansen98,hansen10,mueller}:
If $\mu_j=\mathcal{O}(\rho^{-j})$ with $\rho>1$,
$j=1,2,\ldots,\infty$, then (\ref{eq2}) is severely ill-posed;
if $\mu_j=\mathcal{O}(j^{-\alpha})$, then (\ref{eq2})
is mildly or moderately ill-posed for $\frac{1}{2}<\alpha\le1$ or $\alpha>1$.
Here for mildly ill-posed problems we add the requirement
$\alpha>\frac{1}{2}$, which does not appear in \cite{hofmann86}
but must be met for a linear compact operator equation \cite{hanke93,hansen98}.
In the one dimensional case, i.e., $q=1$, (\ref{eq1})
is severely ill-posed when the kernel function $k(s,t)$ is sufficiently smooth, and
it is moderately ill-posed with $\mu_j=\mathcal{O}(j^{-p-1/2})$,
where $p$ is the highest order of continuous derivatives of
$k(s,t)$; see, e.g., \cite[p.8]{hansen98} and \cite[p.10-11]{hansen10}.
The singular values $\sigma_j$ of discretized problem (\ref{eq1})
resulting from the continuous (\ref{eq2}) inherit
the decay properties of $\mu_j$ \cite{hansen98,hansen10}, provided
that discretizations are fine enough, so that the classification
applies to (\ref{eq1}) as well.

Bj\"{o}rck and Eld\'{e}n in their 1979 survey \cite{bjorck79}
foresightedly expressed a fundamental concern on CGLS (and LSQR): {\em More
research is needed to tell for which problems this approach will work, and
what stopping criterion to choose.} See also \cite[p.145]{hansen98}.
Hanke and Hansen \cite{hanke93} and Hansen \cite{hansen07} address
that a strict proof of the regularizing properties of conjugate gradients is
extremely difficult. Over the years, an enormous effort has been made to the study
of regularizing effects of LSQR and CGLS; see, e.g.,
\cite{firro97,hanke01,hansen98,hansen10,hps16,hps09,huangjia,
jia18a,jia18b,kirsch,nolet,paige06,scales}.
To echo the concern of Bj\"{o}rck
and Eld\'{e}n, such a definition has been introduced in \cite{huangjia,jia18a}:
If a regularized solution to \eqref{eq1} is at least as accurate as
$x_{k_0}^{tsvd}$, then it is called a best possible 2-norm filtering
regularized solution. If the regularized solution by an iterative regularization
solver at semi-convergence is such a best possible one, then
the solver is said to have the {\em full} regularization.
Otherwise, the solver is said to have only
the {\em partial} regularization.

Since it had been unknown whether or not LSQR, CGME and LSMR
have the full regularization for a given \eqref{eq1},
one commonly combines them with some explicit
regularization \cite{aster,hansen98,hansen10}.
The hybrid LSQR variants have been advocated by Bj\"{o}rck and Eld\'{e}n
\cite{bjorck79} and O'Leary and Simmons \cite{oleary81}, and improved and
developed by Bj\"orck \cite{bjorck88}, Bj\"{o}rck, Grimme and
van Dooren \cite{bjorck94}, and Renaut {\em et al}. \cite{renaut}.
A hybrid LSQR first projects \eqref{eq1} onto Krylov
subspaces and then regularizes the projected problems explicitly.
It aims to remove the effects
of small Ritz values and expands Krylov subspaces until they
captures all needed dominant SVD components of $A$
\cite{bjorck88,bjorck94,hanke93,oleary81}, so that
the error norms of regularized
solutions and the residual norms possibly decrease further until they ultimately
stabilize. The hybrid LSQR, CGME and LSMR have been intensively studied in, e.g.,
\cite{berisha,chung08,chung15,hanke01,hanke93,renaut}
and \cite{aster,hansen10}.

If an iterative solver itself, e.g., LSQR,
is theoretically proved and practically identified to
have the full regularization, one simply stops it after
a few iterations of semi-convergence, and no complicated hybrid variant is
needed. In computation, semi-convergence can be in principle determined by
a parameter-choice method, such as the L-curve criterion and the discrepancy
principle. Therefore, we cannot emphasize
too much the importance of proving the full or partial
regularization of LSQR, CGLS, CGME and LSMR. By the definition of
the full or partial regularization, a fundamental question is:
{\em  Do LSQR, CGLS, LSMR and CGME have the full or partial
regularization for severely, moderately and mildly ill-posed problems?}

Regarding LSQR and CGLS, for the three kinds of ill-posed problems
described above, the author \cite{jia18a,jia18b} has proved that
LSQR has the full regularization for severely and moderately ill-posed
problems with certain suitable $\rho>1$ and $\alpha>1$
under the assumption that all the singular values of $A$ are
simple. In applications, there are 2D image belurring problems
where the matrices $A$'s
have multiple singular values \cite{gazzola18}.
In this paper, we extend the results in \cite{jia18a,jia18b} to the
multiple singular value case. In Section~\ref{svdanal}, we
reformulate the TSVD method and standard-form Tikhonov
regularization in the multiple singular value case,
showing that they compute regularized solutions
as if they work on a modified form of \eqref{eq1}, where the coefficient matrix
has the distinct singular values of $A$ as its nonzero simple singular
values. In Section~\ref{lsqr}, we show
that LSQR works as if it solves the same modified one of
\eqref{eq1}. In this way, we build a bridge that connects
the regularizing effects of the TSVD method and those of LSQR, so that we can
analyze the regularization ability of LSQR by taking the best TSVD regularized
solution as the reference standard. In Sections~\ref{sine}--\ref{alphabeta},
we extend the main results in \cite{jia18a,jia18b} to the multiple singular
value case. Consequently, we can draw the same
conclusions as those in \cite{jia18a,jia18b}.

After the above, we consider
some important issues that have received no attention in the literature:
best, near best and general rank $k$ approximations to $A$ for
the ill-posed problems with $\alpha>1$ and $0<\alpha\leq 1$, respectively,
which include mildly ill-posed problems, and some intrinsic relationships
between them and the approximation properties of their nonzero singular values.
These results apply to LSQR, where the Ritz values, i.e.,
the nonzero singular values of rank $k$ approximation
matrices generated by Golub-Kahan bidiagonalization,
critically decide the regularization ability of LSQR.
We will show that, unlike for severely and moderately ill-posed problems
with suitable $\rho>1$ and $\alpha>1$,
a best or near best rank $k$ approximation to $A$ does not mean
that its nonzero singular values approximate the large singular values of $A$
in natural order. Furthermore, for $0<\alpha\leq 1$,
given the accuracy of the rank $k$ approximation in LSQR,
we establish more insightful results on the
nonzero singular values of the rank $k$ approximation matrix,
which estimate their maximum possible number that are smaller than $\sigma_{k+1}$.
These results also apply to the Krylov solvers CGME, MINRES and MR-II,
and GMRES and RRGMRES \cite{hanke95,hansen10}, each of which
generates its own rank $k$ approximation to $A$ at iteration $k$.
All these constitutes the work of Section~\ref{morerank}.
In Section~\ref{exper}, we report the numerical experiments
to confirm the results.
Finally, we conclude the paper in Section~\ref{concl}.

Throughout the paper, we denote by
$\mathcal{K}_{k}(C, w)= span\{w,Cw,\ldots,C^{k-1}w\}$
the $k$ dimensional Krylov subspace generated
by the matrix $\mathit{C}$ and the vector $\mathit{w}$, and by $I$ and
$\mathbf{0}$ the identity matrix
and the zero matrix whose orders are omitted whenever
clear from the context.

\section{The reformulation and analysis of the TSVD method and
standard-form Tikhonov regularization in the multiple singular
value case}\label{svdanal}

In order to extend the results in \cite{jia18a,jia18b} to
the multiple singular value case, we need to
reorganize the SVD of $A$ and reformulate the TSVD method
and standard-form Tikhonov regularization by taking $b$ in \eqref{eq1} into
account carefully.
To this end, we must make numerous necessary changes and preparations,
as will be detailed below.

Let the SVD of $A$ be
\begin{equation}\label{eqsvd}
  A=\widehat{U}\left(\begin{array}{c} \Sigma \\ \mathbf{0} \end{array}\right)
  \widehat{V}^{T},
\end{equation}
where $\widehat{U} = (\widehat{U}_1,\widehat{U}_2,\ldots,\widehat{U}_s,\widehat{U}_{\perp})
\in\mathbb{R}^{m\times m}$ with $\widehat{U}_i\in \mathbb{R}^{m\times c_i}$ and
$\widehat{V} = (\widehat{V}_1,\widehat{V}_2,\ldots,\widehat{V}_s)
\in\mathbb{R}^{n\times n}$
with $\widehat{V}_i\in \mathbb{R}^{n\times c_i}$
are orthogonal, $\Sigma = {\rm diag} (\sigma_1 I_{c_1},\sigma_2 I_{c_2},
\ldots,\sigma_s I_{c_s})$ with the $s$ distinct
singular values $\sigma_1>\sigma_2 >\cdots >\sigma_s>0$, each $\sigma_i$ is $c_i$
multiple and $I_{c_i}$ the $c_i\times c_i$ identify matrix,
and the superscript $T$ denotes
the transpose of a matrix or vector. Then for a given Gaussian white noise $e$, by
\eqref{eqsvd} we obtain
\begin{equation}\label{eq4}
  x_{naive}=\sum\limits_{i=1}^{s}\frac{\widehat{V}_i\widehat{U}_i^{T}b}{\sigma_i} =
  \sum\limits_{i=1}^{s}\frac{\widehat{V}_i\widehat{U}_i^{T}b_{true}}{\sigma_i} +
  \sum\limits_{i=1}^{s}\frac{\widehat{V}_i\widehat{U}_i^{T}e}{\sigma_i}
  =x_{true}+\sum\limits_{i=1}^{s}\frac{\widehat{V}_i\widehat{U}_i^{T}e}{\sigma_i},
\end{equation}
where $x_{true}=A^{\dagger}b_{true}$ with $\|x_{true}\|=\|A^{\dagger}b_{true}\|=
\left(\sum_{i=1}^s\frac{\|\widehat{U}_i^Tb_{true}\|^2}{\sigma_i^2}\right)^{1/2}$,
and the norm of the second term is very large (huge) for $\|e\|$ fixed.

The discrete Picard condition on \eqref{eq1} stems from the
necessary requirement
$$
\|x_{true}\|=\|A^{\dagger}b_{true}\|=
\left(\sum_{i=1}^s\frac{\|\widehat{U}_i^Tb_{true}\|^2}{\sigma_i^2}\right)^{1/2}\leq C
$$
with some constant $C$, independent of $n$ and plays a fundamental role
in the solution of linear discrete ill-posed problems; see, e.g.,
\cite{aster,gazzola15,hansen90b,hansen98,hansen10,kern,renaut}.
It states that, on average, the (generalized) Fourier coefficients
$\|U_i^{T}b_{true}\|$ decay faster than $\sigma_i$, which enables
regularization to compute useful approximations to $x_{true}$.
The following common model has been used throughout Hansen's books
\cite{hansen98,hansen10} and
the references therein as well as \cite{jia18a,jia18b} and the current paper:
\begin{equation}\label{picard}
  \|\widehat{U}_i^T b_{true}\|=\sigma_i^{1+\beta},\ \ \beta>0,\ i=1,2,\ldots,s,
\end{equation}
where $\beta$ is a model parameter that controls the decay rates of
$\|\widehat{U}_i^T b_{true}\|$. We remark that Hansen \cite{hansen98,hansen10}
uses the individual columns of $\widehat{U}_i$ in \eqref{picard} if $\sigma_i$ is
multiple, which is equivalent to assuming that each column of $\widehat{U}_i$ and
the corresponding one of $\widehat{V}_i$ makes the same contributions
to $x_{naive}$ and $x_{true}$ and the TSVD method must
use {\em all} the columns of $\widehat{U}_i$ and $\widehat{V}_i$ associated with $\sigma_i$
to form a regularized solution.
It is trivial to unify the discrete Picard condition on the individual
columns of $\widehat{U}_i$ in the form of \eqref{picard} for a multiple $\sigma_i$.

Based on the above, in the multiple singular value case,
the TSVD method \cite{hansen98,hansen10} solves
\eqref{posed} by dealing with the problem
\begin{equation}\label{tsvd}
\min\|x\| \ \ \mbox{subject to}\ \
 \|A_kx-b\|=\min
\end{equation}
for some $k$, where $A_k$ is a best rank $c_1+c_2+\cdots+c_k$ approximation to $A$
and the most common choice (cf. \cite[p.12]{bjorck96}) is
\begin{equation}\label{rankak}
A_k=(\widehat{U}_1,\widehat{U}_2,\ldots,\widehat{U}_k)\Sigma_k
(\widehat{V}_1,\widehat{V}_2,\ldots,\widehat{V}_k)^T
\end{equation}
with $\Sigma_k={\rm diag} (\sigma_1 I_{c_1},\sigma_2 I_{c_2},
\ldots,\sigma_k I_{c_k})$ and $\|A-A_k\|=\sigma_{k+1}$.

In order to extend the results in \cite{jia18a,jia18b} to
the multiple singular value case, the first key step is to
take the right-hand side $b$ into consideration and to reorganize
\eqref{eqsvd} so as to obtain an SVD of $A$ in some desired form
by selecting a specific set of left and right singular vectors
corresponding to a multiple singular value $\sigma_i$ of $A$.
Specifically,
for the $c_i$ multiple $\sigma_i,\ i=1,2,\ldots,s$, we choose an orthonormal
basis of its left singular subspace by requiring that
$b$ have a nonzero orthogonal projection on just
{\em one} left singular vector $u_i$ in the
singular subspace and {\em no} components in the remaining $c_i-1$ ones.
Precisely, recall that the columns of $\widehat{U}_i$ form an
orthonormal basis of the unique left
singular subspace associated with $\sigma_i$. Then we must have
\begin{equation}\label{newuk}
u_i=\frac{\widehat{U}_i\widehat{U}_i^Tb}{\|\widehat{U}_i^Tb\|},
\end{equation}
where $\widehat{U}_i\widehat{U}_i^T$ is the orthogonal projector onto the
left singular subspace associated with $\sigma_i$. With such $u_i$,
define the corresponding right singular vector
$v_i = A^T u_i/\sigma_i,\ i=1,2,\ldots,s$. We then
select the other $c_i-1$ orthonormal
left singular vectors which are orthogonal to $u_i$ and, together with $u_i$,
form a new orthonormal basis of the left singular subspace associated
with $\sigma_i$. We define the
corresponding $c_i-1$ right singular vectors in the same way
as $v_i$; they, together with $v_i$, form a new orthonormal basis of
the right singular subspace associated with $\sigma_i$.

Write the above resulting new left and right
singular vector matrices as $\widetilde{U}_i$ and $\widetilde{V}_i$
with $u_i$ and $v_i$ as their first columns, respectively.
Then there exist $c_i\times c_i$ orthogonal matrices $Q_{i,l}$
and $Q_{i,r}$ such that $\widetilde{U}_i=\widehat{U}_iQ_{i,l}$ and
$\widetilde{V}_i=\widehat{V}_iQ_{i,r}$.
After such treatment, we obtain a desired compact SVD
\begin{equation}\label{resvd}
A=\widetilde{U}\Sigma \widetilde{V}^T
\end{equation}
with $\widetilde{U}=
(\widetilde{U}_1,\widetilde{U}_2,\ldots,\widetilde{U}_s)$ and
$\widetilde{V}=(\widetilde{V}_1,\widetilde{V}_2,
\ldots,\widetilde{V}_s)$. We remind that
$u_i$ defined above is unique since the orthogonal projection of $b$
onto the left singular subspace associated with $\sigma_i$ is unique
and does not depend on the choice of its orthonormal basis.

Now we need to prove that $\widetilde{U}_i$ satisfies the discrete Picard
condition \eqref{picard}. To see this, notice that
\begin{equation}\label{newpicard}
\|\widetilde{U}_i^T b_{true}\|=\|Q_{i,l}^T\widehat{U}_i^Tb_{true}\|=\|\widehat{U}_i^Tb_{true}\|
=\sigma_i^{1+\beta},\ i=1,2,\ldots,s.
\end{equation}
Particularly, we have
\begin{equation}\label{newpicard2}
|u_i^Tb_{true}|\leq \|\widetilde{U}_i^T b_{true}\|=\sigma_i^{1+\beta},\ i=1,2,\ldots,s.
\end{equation}
We will take the equality when using \eqref{newpicard2} later, which does not
affect all the proofs and results to be presented.

With \eqref{newuk} and  \eqref{resvd}, a crucial observation is
that the solution $x_k^{tsvd}$ to \eqref{tsvd} becomes
\begin{equation}\label{tsvdsolution}
x_k^{tsvd}=A_k^{\dagger}b=\sum_{i=1}^k\frac{u_i^Tb}{\sigma_i}v_i,
\end{equation}
which consists of the first $k$ large {\em distinct} dominant SVD components
$(\sigma_i,u_i,v_i)$ of $A$.

Define the new $m\times n$ matrix
\begin{equation}\label{aprime}
A^{\prime}=U\Sigma^{\prime} V^T,
\end{equation}
where $\Sigma^{\prime}={\rm diag}(\sigma_1,\sigma_2,
\ldots,\sigma_s,\mathbf{0})\in\mathbb{R}^{n\times n}$,
$U=(u_1,u_2,\ldots,u_s, U_{\perp})\in\mathbb{R}^{m\times n}$ and
$V=(v_1,v_2,\ldots,v_s,V_{\perp})\in\mathbb{R}^{n\times n}$
with $U_{\perp}$ and $V_{\perp}$ consisting of the other columns
of $\widetilde{U}$ and $\widetilde{V}$ defined by \eqref{resvd}, respectively.

Write
\begin{equation}\label{uvs}
U=(U_s,U_{\perp}),\ \ V=(V_s,V_{\perp}).
\end{equation}
By construction, $A^{\prime}$ has
the $s$ nonzero simple singular values $\sigma_i,\ i=1,2,\ldots,s$ and
$n-s$ multiple zero singular value, and $U_{\perp}^Tb=\mathbf{0}$.
We then see from \eqref{tsvdsolution} that the best rank $c_1+c_2+\cdots+c_k$
approximation $A_k$ to $A$ in \eqref{tsvd} can be equivalently replaced by the
best rank $k$ approximation
\begin{equation}\label{rankakprime}
A_k^{\prime}=U_k\Sigma_k^{\prime}V_k^T,\ k=1,2,\ldots,s,
\end{equation}
to $A^{\prime}$, where $U_k=(u_1,\ldots,u_k)$,
$V_k=(v_1,\ldots,v_k)$ and $\Sigma_k^{\prime}=
{\rm diag}(\sigma_1,\ldots,\sigma_k)$ because
$x_k^{tsvd}=A_k^{\dagger}b=(A_k^{\prime})^{\dagger}b$.

With the above analysis and simple justifications, we can
present the following theorem.

\begin{theorem}\label{tsvdaa}
Let $A_k$, $A^{\prime}$ and $A_k^{\prime}$ be defined by \eqref{rankak},
\eqref{aprime} and \eqref{rankakprime}.
Then for $k=1,2,\ldots,s$ the TSVD solutions $x_k^{tsvd}$ satisfy
\begin{align}
x_k^{tsvd}&=A_k^{\dagger}b=
(A_k^{\prime})^{\dagger}b,\label{xksolution}\\
\|Ax_k^{tsvd}-b\|&=\|A_kx_k^{tsvd}-b\|=\|A_k^{\prime}x_k^{tsvd}-b\|
=\|A^{\prime}x_k^{tsvd}-b\|. \label{residual}
\end{align}
Particularly, for $k=s$, we have
\begin{equation}\label{naive}
x_{naive}=x_s^{tsvd}=A^{\dagger}b=(A^{\prime})^{\dagger}b.
\end{equation}

The above results show that solving \eqref{tsvd} amounts to solving
the problem
\begin{equation}\label{newtsvd}
\min\|x\| \ \ \mbox{subject to}\ \
 \|A_k^{\prime}x-b\|=\min
\end{equation}
for the same $k=1$ and that \eqref{tsvd} and \eqref{newtsvd}
have the same solutions and residual norms for $k=1,2,\ldots,s$.
\end{theorem}

\begin{remark}
Relations \eqref{xksolution}--\eqref{newtsvd} mean that
the TSVD method for solving \eqref{posed}
works exactly as if it solves the regularization problem
\begin{equation}\label{simposed}
\min\limits_{x\in \mathbb{R}^{n}}\|x\| \ \ \mbox{subject to}\ \
\|A^{\prime}x-b\|\leq \tau \|e\|
\end{equation}
with $\tau\approx 1$, and it computes the same TSVD regularized solutions
$x_k^{tsvd},\ k=1,2,\ldots,s$, to \eqref{eq1} and the modified problem
\begin{equation}\label{simeq1}
\min\limits_{x\in \mathbb{R}^{n}}\|A^{\prime}x-b\|.
\end{equation}
Relation \eqref{newpicard} or \eqref{newpicard2} states that \eqref{simeq1}
satisfies the discrete Picard condition.
\end{remark}

The covariance matrix of the Gaussian white noise $e$
is $\eta^2 I$, and the expected value $\mathcal{E}(\|e\|^2)=m \eta^2$.
With the SVD \eqref{aprime} of $A^{\prime}$, it holds that
$\mathcal{E}(|u_i^Te|)=\eta,\,i=1,2,\ldots,s$, and
$\|e\|\approx \sqrt{m}\eta$ and $|u_i^Te|\approx \eta,\
i=1,2,\ldots,s$; see, e.g., \cite[p.70-1]{hansen98} and \cite[p.41-2]{hansen10}.
The noise $e$ thus affects $u_i^Tb,\ i=1,2,\ldots,s,$ {\em more or less equally}.
Relation \eqref{newpicard2} shows that for large singular
values $|{u_i^{T}b_{true}}|/{\sigma_i}$ is dominant relative to
$|u_i^{T}e|/{\sigma_i}$. Once
$| u_i^T b_{true}| \leq | u_i^T e|$ from some $i$ onwards, the small singular
values magnify $|u_i^{T}e|/{\sigma_i}$, and the noise
$e$ dominates $| u_i^T b|/\sigma_i$ and must be suppressed. The
transition point $k_0$ is such that
\begin{equation}\label{picard1}
| u_{k_0}^T b|\approx | u_{k_0}^T b_{true}|> | u_{k_0}^T e|\approx
\eta, \ | u_{k_0+1}^T b|
\approx | u_{k_0+1}^Te|
\approx \eta;
\end{equation}
see \cite[p.42, 98]{hansen10} and a similar description \cite[p.70-1]{hansen98}.
In this sense, the $\sigma_k$ are divided into the $k_0$ large
and $s-k_0$ small ones. The TSVD solutions
\begin{equation}\label{solution}
  x^{tsvd}_k=\left\{\begin{array}{ll} \sum\limits_{i=1}^{k}\frac{u_i^{T}b}
  {\sigma_i}{v_i}\approx
  \sum\limits_{i=1}^{k}\frac{u_i^{T}b_{true}}
{\sigma_i}{v_i},\ \ \ &k\leq k_0;\\ \sum\limits_{i=1}^{k}\frac{u_i^{T}b}
{\sigma_i}{v_i}\approx
\sum\limits_{i=1}^{k_0}\frac{u_i^{T}b_{true}}{\sigma_i}{v_i}+
\sum\limits_{i=k_0+1}^{k}\frac{u_i^{T}e}{\sigma_i}{v_i},\ \ \ &k_0<k\leq s.
\end{array}\right.
\end{equation}
It is easily justified from \cite[p.70-1]{hansen98} and
\cite[p.71,86-8,96]{hansen10}
that $x_k^{tsvd}$ first converges to $x_{true}$
and the error $\|x_{true}-x_k^{tsvd}\|$ and
the residual norm $\|A^{\prime}x_k^{tsvd}-b\|=\|A_k^{\prime}
x_k^{tsvd}-b\|$ monotonically
decrease until $\|A^{\prime} x_k^{tsvd}-b\|\approx \|e\|$ for $k=k_0$, afterwards
$x_k^{tsvd}$ diverges and instead converges to $x_{naive}$,
while the residual norm $\|A^{\prime} x_k^{tsvd}-b\|=\|A_k^{\prime} x_k^{tsvd}-b\|$
stabilizes
for $k$ not close to $s$. Therefore, the index $k$ plays
the role of the regularization parameter, $x_k^{tsvd}$ exhibits
typical semi-convergence at $k=k_0$, and the best regularized solution
$x_{k_0}^{tsvd}=(A_{k_0}^{\prime})^{\dagger}b$ has minimum 2-norm error.

By the construction of $U$ and $V$, it follows from \eqref{eqsvd} and
\eqref {aprime} that, for a given
parameter $\lambda$, the solution $x_{\lambda}$ of the Tikhonov regularization
\eqref{tikhonov} is
\begin{equation}\label{eqfilter}
  x_{\lambda} = \sum\limits_{i=1}^{s}f_i\frac{\widehat{V}_i\widehat{U}_i^{T}b}
  {\sigma_i}= \sum\limits_{i=1}^{s}f_i\frac{u_i^{T}b}{\sigma_i}v_i,
\end{equation}
which is a filtered SVD expansion of $A^{\prime}$,
where $f_i=\frac{\sigma_i^2}{\sigma_i^2+\lambda^2},\ i=1,2,\ldots,s$, are
called filters. The above relation has proved the following result.

\begin{theorem}\label{simtikh}
For the same $\lambda>0$, \eqref{tikhonov} with $L=I$ is equivalent
to the standard-form Tikhonov regularization
\begin{equation}\label{mtikhonov}
  \min\limits_{x\in \mathbb{R}^{n}}\{\|A^{\prime}x-b\|^2+\lambda^2\|x\|^2\}
\end{equation}
of the modified \eqref{simeq1}.
\end{theorem}

$x_k^{tsvd}$ by the TSVD method is a special filtered SVD expansion,
where $f_i=1,\ i=1,2,\ldots,k$ and
$f_i=0,\ i=k+1,\ldots,s$. The best Tikhonov regularized
solution $x_{\lambda_{opt}}$, which is defined as
$\|x_{true}-x_{\lambda_{opt}}\|=\min_{\lambda\geq 0}\|x_{true}-x_{\lambda}\|$,
retains the $k_0$ dominant SVD components
of $A^{\prime}$ and dampens the other $s-k_0$ small SVD components as much as
possible. The semi-convergence of the Tikhonov regularization
method occurs at $\lambda_{opt}$ when the parameter $\lambda$ varies from zero
to infinity.

Finally, we stress that our above changes and reformulations are purely for a
mathematical analysis, which aims
to extend the results in \cite{jia18a,jia18b}
to the multiple singular value case. Computationally, we
never need to reorganize the SVD of $A$ and
construct $A^{\prime}$.

\section{The LSQR algorithm}\label{lsqr}

The LSQR algorithm is based on Golub-Kahan bidiagonalization,
Algorithm~\ref{alg:lb}, that
computes two orthonormal bases $\{q_1,q_2,\dots,q_k\}$ and
$\{p_1,p_2,\dots,p_{k+1}\}$  of $\mathcal{K}_{k}(A^{T}A,A^{T}b)$ and
$\mathcal{K}_{k+1}(A A^{T},b)$  for $k=1,2,\ldots,n$,
respectively.
\begin{algorithm}
\begin{itemize}
   \item[1.] Take $ p_1=b/\|b\| \in \mathbb{R}^{m}$, and define
   $\beta_1{q_0}=\mathbf{0}$ with $\beta_1=\|b\|$.

   \item[2.] For $j=1,2,\ldots,k$
\begin{description}
  \item[(i)]
  $r = A^{T}p_j - \beta_j{q_{j-1}}$
  \item[(ii)] $\alpha_j = \|r\|;q_j = r/\alpha_j$
  \item[(iii)]
   $   z = Aq_j - \alpha_j{p_{j}}$
  \item[(iv)]
  $\beta_{j+1} = \|z\|;p_{j+1} = z/\beta_{j+1}.$
\end{description}
\end{itemize}
\label{alg:lb}
\end{algorithm}

Algorithm~\ref{alg:lb} can be written in the matrix form
\begin{align}
  AQ_k&=P_{k+1}B_k,\label{eqmform1}\\
  A^{T}P_{k+1}&=Q_{k}B_k^T+\alpha_{k+1}q_{k+1}(e_{k+1}^{(k+1)})^{T},\label{eqmform2}
\end{align}
where $e_{k+1}^{(k+1)}$ is the $(k+1)$-th canonical basis vector of
$\mathbb{R}^{k+1}$, $P_{k+1}=(p_1,p_2,\ldots,p_{k+1})$,
$Q_k=(q_1,q_2,\ldots,q_k)$, and
\begin{equation}\label{bk}
  B_k = \left(\begin{array}{cccc} \alpha_1 & & &\\ \beta_2 & \alpha_2 & &\\ &
  \beta_3 &\ddots & \\& & \ddots & \alpha_{k} \\ & & & \beta_{k+1}
  \end{array}\right)\in \mathbb{R}^{(k+1)\times k}.
\end{equation}
It follows from \eqref{eqmform1} that
\begin{equation}\label{Bk}
B_k=P_{k+1}^TAQ_k.
\end{equation}
We remind that the singular values $\theta_i^{(k)},\,i=1,2,\ldots,k$
of $B_k$, called the Ritz values of $A$ with
respect to the left and right subspaces $span\{P_{k+1}\}$ and $span\{Q_k\}$,
are all simple, provided that Algorithm~\ref{alg:lb} does not break down
until step $k$.

At iteration $k$, LSQR solves the problem
$$
\|Ax_k^{lsqr}-b\|=\min_{x\in \mathcal{K}_k(A^TA,A^Tb)}
\|Ax-b\|
$$
and computes the iterate $x_k^{lsqr}=Q_ky_k^{lsqr}$ with
\begin{equation}\label{yk}
  y_k^{lsqr}=\arg\min\limits_{y\in \mathbb{R}^{k}}\|B_ky-\beta_1 e_1^{(k+1)}\|
  =\beta_1  B_k^{\dagger} e_1^{(k+1)},
\end{equation}
where $e_1^{(k+1)}$ is the first canonical basis vector of $\mathbb{R}^{k+1}$,
and the solution norm $\|x_k^{lsqr}\|=\|y_k^{lsqr}\|$ increases
and the residual norm $\|Ax_k^{lsqr}-b\|=\|B_ky_k^{lsqr}-\beta_1 e_1^{(k+1)}\|$
decreases monotonically with respect to
$k$. From $\beta_1 e_1^{(k+1)}=P_{k+1}^T b$ and \eqref{yk}, we obtain
\begin{equation}\label{xk}
x_k^{lsqr}=Q_k B_k^{\dagger} P_{k+1}^Tb,
\end{equation}
which solves the problem
$$
\min\|x\| \ \ \mbox{ subject to }\ \ \|P_{k+1}B_kQ_k^Tx-b\|=\min.
$$

Next we will prove that Algorithm~\ref{alg:lb} and LSQR work
exactly as if they are applied to $A^{\prime}$ and  \eqref{simeq1},
that is, they generate the same results when applied to solving Problems
\eqref{eq1} and \eqref{simeq1}.

By \eqref{uvs}, let us expand $b$ as
$$
b=U_sU_s^Tb+(I-U_sU_s^T)b=\sum_{j=1}^s\xi_j u_j+(I-U_sU_s^T)b.
$$
By the SVD \eqref{resvd} of $A$ and the SVD \eqref{aprime} of $A^{\prime}$
as well as the description on them, it is straightforward to justify that
\begin{equation}\label{ata}
\mathcal{K}_k(A^TA, A^Tb)=\mathcal{K}_k((A^{\prime})^TA^{\prime},(A^{\prime})^Tb)
\end{equation}
and
\begin{equation}\label{aat}
\mathcal{K}_k(AA^T, b)=\mathcal{K}_k(A^{\prime}(A^{\prime})^T,b)
\end{equation}
by noting that
\begin{equation}\label{aprimea}
(A^TA)^i A^Tb=\left((A^{\prime})^TA^{\prime}\right)^i(A^{\prime})^Tb
=\sum_{j=1}^s \xi_j \sigma_j^{2i+1}v_j
\end{equation}
for any integer $i\geq 0$ and
\begin{equation}\label{aaprime}
(AA^T)^i b=\left(A^{\prime}(A^{\prime})^T\right)^ib
=\sum_{j=1}^s \xi_j \sigma_j^{2i}u_j
\end{equation}
for any integer $i\geq 1$. Thus, for the given $b$, Algorithm~\ref{alg:lb}
works on $A$ exactly as if it does on $A^{\prime}$, that is,
\eqref{eqmform1}--\eqref{Bk} hold when $A$ is replaced by $A^{\prime}$. As
a result, the $k$ distinct Ritz values $\theta_i^{(k)}$
approximate $k$ nonzero singular values of $A^{\prime}$, i.e., $k$ distinct
singular values of $A$. Particularly, from \eqref{Bk} we have
\begin{equation}\label{newbk}
B_k=P_{k+1}^TA^{\prime}Q_k.
\end{equation}
Moreover, \eqref{aprimea} and \eqref{aaprime} show
$$
\mathcal{K}_{s+1}((A^{\prime})^TA^{\prime},
(A^{\prime})^Tb)=\mathcal{K}_s((A^{\prime})^TA^{\prime},
(A^{\prime})^Tb),\
\mathcal{K}_{s+2}(A^{\prime}(A^{\prime})^T,b)=
\mathcal{K}_{s+1}(A^{\prime}(A^{\prime})^T,b).
$$
As a result, since $(A^{\prime})^Tb$ has nonzero components in all the
right singular vectors $v_1,v_2,\ldots,v_s$ of $A^{\prime}$ associated
with its nonzero distinct singular values $\sigma_1,\sigma_2,\ldots,
\sigma_s$, Golub-Kahan bidiagonalization cannot break down until step $s$,
and the $s$ singular values $\theta_i^{(s)}$ of $B_s$ are exactly the singular
values $\sigma_1,\sigma_2,\ldots,\sigma_s$ of $A^{\prime}$.
At step $s$, Golub-Kahan bidiagonalization generates the orthonormal
$P_{s+1}$,$Q_s$ and the matrix\footnote{
If $m=n$, it is easily justified that $\beta_{s+1}=0$,
Algorithm~\ref{alg:lb} produces the orthonormal matrices $P_s$, $Q_s$ and
the $s\times s$ lower bidiagonal $B_s$ with the positive diagonals
$\alpha_i$ and subdiagonals $\beta_i$. This does not affect
all the derivation and results followed, and we only need to replace
$P_{s+1}$ by $P_s$. For example, $span\{U_s\}=span\{P_s\}$ in \eqref{uv}.}
\begin{align}
P_{s+1}^T AQ_s&=P_{s+1}^T A^{\prime} Q_s=B_s \label{lbaprime}
\end{align}
and
\begin{align}
span\{V_s\}&=span\{Q_s\},\ span\{U_s\}\subset
span\{P_{s+1}\}. \label{uv}
\end{align}

Since \eqref{eqmform1}--\eqref{Bk}
hold when $A$ is replaced by $A^{\prime}$, just as the TSVD
method (cf. \eqref{newtsvd}), LSQR works exactly as if it
solves \eqref{simposed}. We summarize the results as follows.

\begin{theorem}\label{lsqraa}
The LSQR iterate $x_k^{lsqr}$ is the solution
to  the problem
\begin{equation}\label{lsqrreg}
\min\|x\| \ \ \mbox{ subject to }\ \ \|P_{k+1}B_kQ_k^Tx-b\|=\min
\end{equation}
starting with $k=1$ onwards, and it is a regularized solution
of \eqref{simeq1} and thus of \eqref{eq1} and satisfies
\begin{equation}\label{resaa}
\|Ax_k^{lsqr}-b\|=\|A^{\prime}x_k^{lsqr}-b\|=\|B_ky_k^{lsqr}-\beta_1 e_1^{(k+1)}\|,\
k=1,2,\ldots,s
\end{equation}
with $y_k^{lsqr}$ defined by \eqref{yk}.
\end{theorem}

The rank $k$ approximation $P_{k+1}B_kQ_k^T$ to $A^{\prime}$
in \eqref{lsqrreg} plays a role similar to the best rank $k$ approximation
$A_k^{\prime}$ to $A^{\prime}$ in \eqref{newtsvd}. Recall that the
best rank $k$ approximation $A_k^{\prime}$ to $A^{\prime}$ satisfies
$\|A^{\prime}-A_k^{\prime}\|=\sigma_{k+1}$. As a result,
if $P_{k+1}B_k Q_k^T$ is a near best rank $k$ approximation
to $A^{\prime}$ with an approximate accuracy $\sigma_{k+1}$ and
the $k$ singular values of $B_k$ approximate the first $k$
large ones of $A^{\prime}$ in natural order for $k=1,2,\ldots,k_0$,
that is, they interlace the first $k+1$ large $\sigma_i$
for $k=1,2,\ldots,k_0$,
LSQR has the same regularization ability as the TSVD method
and has the full regularization
because (i) $x_k^{tsvd}$ and $x_k^{lsqr}$ are the regularized
solutions to the two perturbed problems of \eqref{simeq1} that replace
$A^{\prime}$ by the two rank $k$ approximations with the same quality
to $A^{\prime}$, respectively;
(ii) $x_k^{tsvd}$ and $x_k^{lsqr}$ solve the two
essentially same regularization
problems \eqref{newtsvd} and \eqref{lsqrreg}, respectively.
Therefore, the near best rank $k$ approximation of $P_{k+1}B_k Q_k^T$
to $A^{\prime}$ and the approximations of
the $k$ singular values of $B_k$ to the
large ones of $A^{\prime}$ in natural order for $k=1,2,\ldots, k_0$ are
{\em sufficient} conditions for LSQR to have the full regularization.
We will give the precise definition of a near best rank $k$ approximation
to $A^{\prime}$ later.
However, one must be well aware that they are {\em not
necessary} conditions for the full regularization of LSQR, as has
been addressed in \cite{jia18a,jia18b}.

\section{$\sin\Theta$ theorems for the distances between
$\mathcal{K}_k((A^{\prime})^TA^{\prime},(A^{\prime})^Tb)$ and
the dominant right singular subspace $span\{V_k\}$}
\label{sine}

In the multiple singular value case, based on
the work in Sections~\ref{svdanal}--\ref{lsqr}, just as \cite{jia18a,jia18b},
under the discrete Picard condition \eqref{newpicard2},
a complete understanding of the regularization
of LSQR includes accurate solutions of the following problems: (i)
How accurately does $\mathcal{K}_{k}((A^{\prime})^TA^{\prime},(A^{\prime})^Tb)$
approximate the $k$ dimensional dominant right singular subspace $span\{V_k\}$ of
$A^{\prime}$ spanned by the columns of $V_k=(v_1,v_2,\ldots,v_k)$?
(ii) How accurate is the rank $k$ approximation $P_{k+1}B_kQ_k^T$ to
$A^{\prime}$? (iii) When do the $k$ Ritz values $\theta_i^{(k)}$
approximate the the first $k$ large $\sigma_i$ in natural order?
(iv) When does at least a small Ritz value appear, i.e.,
$\theta_k^{(k)}<\sigma_{k+1}$ for some $k\leq k^*$
with $k^*$ the iteration at which the semi-convergence of LSQR
occurs? (v) Does LSQR have the full or partial regularization
when the $k$ Ritz values $\theta_i^{(k)}$ do not approximate the large
singular values of $A$ in natural order for some $k\leq k^*$?

We will focus on Problems (i)-(iv) and extend all the results
in \cite{jia18a,jia18b} to the multiple singular value case.
On the other hand, as one of the main contributions in
this paper, we will make a novel general analysis that covers
but is not limited to Problem (iii)-(iv) and
get more insight into them when $A$ has simple or multiple
singular values.

Based on a well-known result (cf. e.g., van der Sluis and
van der Vorst \cite[Property 2.8]{vorst86}), it is
straightforward to establish the following result,
which, based on the work of Section~\ref{lsqr},
holds when $A$ is replaced by $A^{\prime}$, and
has been used in Hansen \cite{hansen98} and the references therein as well
as in \cite{jia18a}
to illustrate the regularizing effects of LSQR.

\begin{proposition}\label{help}
LSQR with the starting vector $p_1=b/\|b\|$ and CGLS
applied to the normal equation $(A^{\prime})^TA^{\prime}x=(A^{\prime})^Tb$
of \eqref{simeq1} with the zero starting vector
generate the same iterates
\begin{equation}\label{eqfilter2}
  x_k^{lsqr}=\sum\limits_{i=1}^sf_i^{(k)}\frac{u_i^{T}b}{\sigma_i}v_i,\
  k=1,2,\ldots,s,
\end{equation}
where the filters
\begin{equation}\label{filter}
f_i^{(k)}=1-\prod\limits_{j=1}^k\frac{(\theta_j^{(k)})^2-\sigma_i^2}
{(\theta_j^{(k)})^2},\ i=1,2,\ldots,s,
\end{equation}
and the $\theta_j^{(k)}$ are the singular values of $B_k$
labeled as $\theta_1^{(k)}>\theta_2^{(k)}>\cdots>\theta_k^{(k)}$.
\end{proposition}

Relation \eqref{eqfilter2} shows that $x_k^{lsqr}$ has a filtered SVD expansion
similar to \eqref{eqfilter}. It is easily justified that if all
the Ritz values $\theta_j^{(k)}$ approximate the first $k$ singular values
$\sigma_j$ of $A^{\prime}$ in natural order then $f_i^{(k)}\approx 1,\,
i=1,2,\ldots,k$ and the other $f_i^{(k)}$ monotonically approach zero
for $i=k+1,\ldots,s$. This indicates that if the $\theta_j^{(k)}$
approximate the first $k$ singular values
$\sigma_j$ of $A^{\prime}$ in natural order for $k=1,2,\ldots,k_0$ then
$x_{k_0}^{lsqr}$ is accurate as $x_{k_0}^{tsvd}$, meaning that
LSQR has the full regularization and computes
a best possible 2-norm filtering regularized solution.
Using the same proof as that
of \cite[Theorem 3.1]{jia18a}, we obtain the following basic results.

\begin{theorem}\label{propaprime}
The semi-convergence of LSQR must occur at some iteration
$$
k^*\leq k_0.
$$
If the Ritz values $\theta_i^{(k)}$
do not converge to the large singular values of
$A^{\prime}$ in natural order for some $k\leq k^*$, then
$k^*<k_0$ strictly.  On the other hand,
if $k^*<k_0$, then the Ritz values $\theta_j^{(k)}$ must
not converge to the first
$k$ large singular values $\sigma_j$ of $A$ in natural order for some
$k\leq k^*$.
\end{theorem}


The approximation accuracy of $P_{k+1}B_kQ_k^T$ to $A^{\prime}$
and the approximation properties of $\theta_i^{(k)},\ i=1,2,\ldots,k$,
critically depend on how the underlying
$k$ dimensional $\mathcal{K}_{k}((A^{\prime})^TA^{\prime},(A^{\prime})^Tb)$,
from which the iterate $x_k^{lsqr}$ is extracted,
approximates the $k$ dimensional dominant right singular subspace $span\{V_k\}$
of $A^{\prime}$.
In terms of the canonical angles $\Theta(\mathcal{X},\mathcal{Y})$ between
two subspaces $\mathcal{X}$ and $\mathcal{Y}$ of equal
dimension (cf. \cite[p.74-5]{stewart01} and \cite[p.43]{stewartsun}), we present
the following general result, which is the same as Lemma 4.1
in \cite{jia18a} except that $n$ is replaced by $s$.

\begin{lemma}\label{lemma1}
For $k=1,2,\ldots,s-1$ we have
\begin{align}
\|\sin\Theta(\mathcal{V}_k,\mathcal{V}_k^R)\|&=
\frac{\|\Delta_k\|}{\sqrt{1+\|\Delta_k\|^2}}
\label{deltabound}
\end{align}
with $\Delta_k \in \mathbb{R}^{(n-k)\times k}$ defined by \eqref{defdelta},
i.e.,
\begin{equation}\label{tangent}
\|\tan\Theta(\mathcal{V}_k,\mathcal{V}_k^R)\|=
\|\Delta_k\|.
\end{equation}
\end{lemma}

{\em Proof}. In order to prove \eqref{deltabound}, we need to make some
nontrivial changes in the proof of the same result in \cite{jia18a}.

From \eqref{aprime}, observe the Krylov subspace
$\mathcal{K}_{k}((\Sigma^{\prime})^2,
\Sigma^{\prime} U^Tb)=span\{\hat{D}\hat{T}_k\}$ with
$$
  \hat{D}={\rm diag}(\sigma_1 u_1^Tb,\ldots,\sigma_s u_s^Tb,
  \mathbf{0})=\left(\begin{array}{cc}
  D& \\
  & \mathbf{0}
  \end{array}\right)\in\mathbb{R}^{n\times n}
$$
and
$$
  \hat{T}_k=\left(\begin{array}{cccc} 1 &
  \sigma_1^2&\ldots & \sigma_1^{2k-2}\\
1 &\sigma_2^2 &\ldots &\sigma_2^{2k-2} \\
\vdots & \vdots&&\vdots\\
1 &\sigma_s^2 &\ldots &\sigma_s^{2k-2}\\
0 & 0 &\ldots& 0\\
\vdots & \vdots&&\vdots\\
0 & 0 &\ldots& 0
\end{array}\right)=\left(\begin{array}{c}
T_k\\
\mathbf{0}
\end{array}\right)\in \mathbb{R}^{n\times k}.
$$
Partition the diagonal matrix $D$ and the matrix $T_k$ as
\begin{equation*}
  D=\left(\begin{array}{cc} D_1 & 0 \\ 0 & D_2 \end{array}\right)
  \in \mathbb{R}^{s\times s},\ \ \
  T_k=\left(\begin{array}{c} T_{k1} \\ T_{k2} \end{array}\right)
  \in \mathbb{R}^{s\times k},
\end{equation*}
where $D_1,\ T_{k1}\in\mathbb{R}^{k\times k}$. Since $T_{k1}$ is
a Vandermonde matrix with $\sigma_j$ distinct for $j=1,2,\ldots,k$, it is
nonsingular. Therefore, from
$$
\mathcal{V}_k^R=\mathcal{K}_k((A^{\prime})^TA^{\prime},
(A^{\prime})^Tb)=span\{V\hat{D}\hat{T}_k\}
$$
and the structures of $\hat{D}$ and $\hat{T}_k$ as well as \eqref{uvs},
we obtain
\begin{equation*}
\mathcal{V}_k^R=span\{V_sD T_k\}=span \left\{V_s\left
(\begin{array}{c} D_1T_{k1} \\ D_2T_{k2} \end{array}\right)\right\}
=span\left\{V_s\left(\begin{array}{c} I \\ \Delta_k \end{array}\right)\right\},
\end{equation*}
with
\begin{equation}\label{defdelta}
\Delta_k=D_2T_{k2}T_{k1}^{-1}D_1^{-1},
\end{equation}
meaning that $\mathcal{V}_k^R$ is orthogonal to
$V_{\perp}$ in \eqref{uvs}.

Recall \eqref{aprime} and write
\begin{equation}\label{parti}
V_s=(V_k, V_k^{\perp}).
\end{equation}
Define
\begin{equation*}
Z_k=V_s\left(\begin{array}{c} I \\ \Delta_k \end{array}\right)
=V_k+V_k^{\perp}\Delta_k.
\end{equation*}
Then $Z_k^TZ_k=I+\Delta_k^T\Delta_k$, and the
columns of $\hat{Z}_k=Z_k(Z_k^TZ_k)^{-\frac{1}{2}}$
form an orthonormal basis of $\mathcal{V}_k^R$. Therefore,
we have $V_{\perp}^T\hat{Z}_k=\mathbf{0}$ and obtain an orthogonal
direct sum decomposition
\begin{equation*}
\hat{Z}_k=(V_k+V_k^{\perp}\Delta_k)(I+\Delta_k^T\Delta_k)^{-\frac{1}{2}}.
\end{equation*}

Based on the above, \eqref{uvs} and \eqref{parti}, by
the definition of  $\|\sin\Theta(\mathcal{V}_k, \mathcal{V}_k^R)\|$
and $V_{\perp}^T\hat{Z}_k=\mathbf{0}$, we obtain
\begin{align*}
   \|\sin\Theta(\mathcal{V}_k,\mathcal{V}_k^R)\|
   &=\|(V_k^{\perp},V_{\perp})^T\hat{Z}_k\|
   =\|\Delta_k(I+\Delta_k^T\Delta_k)^{-\frac{1}{2}}\|
   =\frac{\|\Delta_k\|}{\sqrt{1+\|\Delta_k\|^2}},
\end{align*}
which proves \eqref{deltabound}. Relation \eqref{tangent} follows
from \eqref{deltabound} directly.
\qquad\endproof

The following theorem gives accurate estimates for $\|\Delta_k\|$ for severely
ill-posed problems.

\begin{theorem}\label{thm2}
Let the SVD of $A^{\prime}$ be as \eqref{aprime}, and
denote $\mathcal{V}_k=span\{V_k\}$
and $\mathcal{V}_k^R=\mathcal{K}_{k}((A^{\prime})^TA^{\prime},(A^{\prime})^Tb)$,
and assume that \eqref{eq1} is
severely ill-posed with $\sigma_j=\mathcal{O}(\rho^{-j})$ and $\rho>1$,
$j=1,2,\ldots,s$, and the discrete Picard condition \eqref{newpicard2} is satisfied.
Then
\begin{align}
\|\Delta_1\|&\leq \frac{\sigma_{2}}{\sigma_1}\frac{\min_{2\leq j\leq n}|u_i^Tb|}{|u_1^Tb|}
\left(1+\mathcal{O}(\rho^{-2})\right),\label{k1}\\
  \|\Delta_k\|&\leq
  \frac{\sigma_{k+1}}{\sigma_k}\frac{\max_{k+1\leq j\leq n}|u_i^Tb|}{\min_{1\leq i\leq k}|u_i^Tb|}
  \left(1+\mathcal{O}(\rho^{-2})\right)
  |L_{k_1}^{(k)}(0)|,\ k=2,3,\ldots,s-1, \label{eqres1}
\end{align}
where
\begin{equation}\label{lk}
|L_{k_1}^{(k)}(0)|=\max_{j=1,2,\ldots,k}|L_j^{(k)}(0)|,
\ |L_j^{(k)}(0)|=\prod\limits_{i=1,i\ne j}^k\frac{\sigma_i^2}{|\sigma_j^2-
\sigma_i^2|},\,j=1,2,\ldots,k.
\end{equation}
\end{theorem}

{\em Proof}. The proofs of the results follow the corresponding ones of
Theorem 4.2 in \cite{jia18a} step by step and are thus omitted.
\qquad\endproof

Relation \eqref{deltabound} is independent
of the degree of ill-posedness of problem \eqref{eq1}.
The following accurate estimates for $|L_{k_1}^{(k)}(0)|$ and
all $|L_j^{(k)}(0)|$, $j=1,2,\ldots,k$ defined by \eqref{lk}
are straightforward from Theorem 4.3 of \cite{jia18a}.

\begin{theorem}\label{estlk}
For the severely ill-posed problem with the singular values
$\sigma_j=\mathcal{O}(\rho^{-j})$ and
suitable $\rho>1$, $j=1,2,\ldots,n$ and $k=2,3,\ldots,s-1$, we have
\begin{align}
|L_k^{(k)}(0)|&=1+\mathcal{O}(\rho^{-2}), \label{lkkest}\\
|L_j^{(k)}(0)|&=\frac{1+\mathcal{O}(\rho^{-2})}
{\prod\limits_{i=j+1}^k\left(\frac{\sigma_{j}}{\sigma_i}\right)^2}
=\frac{1+\mathcal{O}(\rho^{-2})}{\mathcal{O}(\rho^{(k-j)(k-j+1)})},
\ j=1,2,\ldots,k-1, \label{lj0}\\
|L_{k_1}^{(k)}(0)|&=\max_{j=1,2,\ldots,k}|L_j^{(k)}(0)|
=1+\mathcal{O}(\rho^{-2}). \label{lkk}
\end{align}
\end{theorem}

For moderately and mildly ill-posed
problems, the estimates for $\|\Delta_k\|,\ k=1,2,\ldots,s-1$
and the proofs are the same as those of Theorem 4.4 in \cite{jia18a}
except that $n$ is replaced by $s$.

\begin{theorem}\label{moderate}
Assume that \eqref{eq1} is moderately or mildly ill-posed with $\sigma_j=
\zeta j^{-\alpha},\ j=1,2,\ldots,s$,
where $\alpha>\frac{1}{2}$ and $\zeta>0$ is some constant,
and the other assumptions and notation are the same as in Theorem~\ref{thm2}.
Then \eqref{deltabound} holds with
\begin{align}
\|\Delta_1\|&\leq \frac{\min_{2\leq i\leq n}|u_i^Tb|}{| u_1^Tb|}
\sqrt{\frac{1}{2\alpha-1}},\label{mod1}\\
\|\Delta_k\|&\leq \frac{\max_{k+1\leq j\leq n}|u_i^Tb|}{\min_{1\leq i\leq k}|u_i^Tb|}
\sqrt{\frac{k^2}{4\alpha^2-1}+\frac{k}{2\alpha-1}}
|L_{k_1}^{(k)}(0)|,\ k=2,3,\ldots,s-1. \label{mod2}
\end{align}
\end{theorem}

In Theorem 4.5 of \cite{jia18a}, the author has given estimates for
$|L_j^{(k)}(0)|$,
$j=1,2,\ldots,k$ and $|L_{k_1}^{(k)}(0)|$, which carry over to the
multiple singular value case trivially.

\begin{theorem}\label{estlk2}
For the moderately and mildly ill-posed problems with
$\sigma_i=\zeta i^{-\alpha},\ i=1,2,\ldots,n$ and suitable
$\alpha>1$, we have
\begin{align}
|L_j^{(k)}(0) |&\approx\left(1+\frac{j}{2\alpha+1}\right)
\prod_{i=j+1}^k\left(\frac{j}{i}\right)^{2\alpha},\,j=1,2,\ldots,k-1,
\label{lkjmod1}\\
\frac{k}{2\alpha+1}&<|L_{k_1}^{(k)}(0)|\approx 1+\frac{k}{2\alpha+1}
\label{lk1size}
\end{align}
with the lower bound requiring $k$ satisfying $\frac{2\alpha+1}{k}\leq 1$;
for $\frac{1}{2}<\alpha\leq 1$ and $k$ satisfying $\frac{2\alpha+1}{k}\leq 1$,
we have
\begin{equation}\label{lk1sizemild}
\frac{k}{2\alpha+1}<|L_{k_1}^{(k)}(0)|.
\end{equation}
\end{theorem}

The author in \cite{jia18a} has
investigated how $\|\sin\Theta(\mathcal{V}_k,\mathcal{V}_k^R)\|$
affects the smallest Ritz value $\theta_k^{(k)}$ in the simple singular
value case. We can extend the results to the multiple singular
value case in the same form by modifying the proof.

\begin{theorem}\label{initial}
Let $\|\sin\Theta(\mathcal{V}_k,\mathcal{V}_k^R)\|^2=1-\varepsilon_k^2$ with
$0< \varepsilon_k< 1$, $k=1,2,\ldots,s-1$, and
$\tilde{q}_k\in\mathcal{V}_k^R$ with $\|\tilde{q}_k\|=1$ be the vector
having the smallest angle with $span\{(V_k^{\perp},V_{\perp})\}$ defined by
\eqref{uvs} and \eqref{parti}, i.e., the orthogonal complement of $\mathcal{V}_k$
with respect to $\mathbb{R}^n$. Then it holds that
\begin{equation}\label{rqi}
\varepsilon_k^2\sigma_k^2+
(1-\varepsilon_k^2)\sigma_n^2< \tilde{q}_k^T(A^{\prime})^TA^{\prime}\tilde{q}_k<
\varepsilon_k^2\sigma_{k+1}^2+
(1-\varepsilon_k^2)\sigma_1^2.
\end{equation}
If $\varepsilon_k\geq \frac{\sigma_{k+1}}{\sigma_k}$,
then
\begin{equation}
\sqrt{\tilde{q}_k^T(A^{\prime})^TA^{\prime}\tilde{q}_k}>\sigma_{k+1};
\label{est1}
\end{equation}
if $\varepsilon_k^2\leq\frac{\delta}
{(\frac{\sigma_1}{\sigma_{k+1}})^2-1}$ for a given arbitrarily small
$\delta>0$, then
\begin{equation}\label{thetasigma}
\theta_k^{(k)}<(1+\delta)^{1/2}\sigma_{k+1},
\end{equation}
meaning that $\theta_k^{(k)}<\sigma_{k+1}$
once $\varepsilon_k$ is sufficiently small, i.e.,
$\|\sin\Theta(\mathcal{V}_k,\mathcal{V}_k^R)\|$ is sufficiently close to
one.
\end{theorem}

{\em Proof}.
Since the columns of $Q_k$ generated by Golub-Kahan bidiagonalization form an
orthonormal basis of $\mathcal{V}_k^R$, by definition and the assumption on
$\tilde{q}_k$ we have
\begin{align}
\|\sin\Theta(\mathcal{V}_k,\mathcal{V}_k^R)\|&=\|(V_k^{\perp},V_{\perp})^TQ_k\|
=\|(V_k^{\perp},V_{\perp})(V_k^{\perp},V_{\perp})^TQ_k\| \notag\\
&=\max_{\|c\|=1}\|(V_k^{\perp},V_{\perp})(V_k^{\perp},V_{\perp})^TQ_kc\|\notag\\
&=\|(V_k^{\perp},V_{\perp})(V_k^{\perp},V_{\perp})^T Q_kc_k\| \notag\\
&=\|(V_k^{\perp},V_{\perp})(V_k^{\perp},V_{\perp})^T\tilde{q}_k\|
=\|(V_k^{\perp},V_{\perp})^T\tilde{q}_k\|
=\sqrt{1-\varepsilon_k^2}
\label{qktilde}
\end{align}
with $\tilde{q}_k=Q_kc_k\in\mathcal{V}_k^R$ and $\|c_k\|=1$.

Notice that $V=(V_k,V_k^{\perp},V_{\perp})$.
Expand $\tilde{q}_k$ as the following orthogonal direct sum decomposition:
\begin{equation}\label{decompqk}
\tilde{q}_k=(V_k^{\perp},V_{\perp})(V_k^{\perp},V_{\perp})^T\tilde{q}_k
+V_kV_k^T\tilde{q}_k.
\end{equation}
Then from $\|\tilde{q}_k\|=1$ and \eqref{qktilde} we obtain
\begin{align}\label{angle2}
\|V_k^T\tilde{q}_k\|&=\|V_kV_k^T\tilde{q}_k\|=
\sqrt{1-\|(V_k^{\perp},V_{\perp})(V_k^{\perp},
V_{\perp})^T\tilde{q}_k\|^2}=\sqrt{1-(1-\varepsilon_k^2)}=\varepsilon_k.
\end{align}

We next bound the Rayleigh quotient of $(A^{\prime})^TA^{\prime}$
with respect to $\tilde{q}_k$ from below. By
$A^{\prime}=U\Sigma^{\prime} V^T$ defined in
\eqref{aprime} and  \eqref{parti},
we partition
$$
\Sigma^{\prime}=\left(\begin{array}{ccc}
\Sigma_k^{\prime} & &\\
&\Sigma_{k,\perp}^{\prime}&\\
& & \mathbf{0}
\end{array}
\right),
$$
where $\Sigma_k^{\prime}={\rm diag}(\sigma_1,\sigma_2,\ldots,\sigma_k)$ and
$\Sigma_{k,\perp}^{\prime}={\rm diag}(\sigma_{k+1},\sigma_{k+2},\ldots,\sigma_s)$.
Making use of $(A^{\prime})^TA^{\prime}V_k=V_k(\Sigma_k^{\prime})^2$,
$(A^{\prime})^TA^{\prime}V_k^{\perp}=
V_k^{\perp}(\Sigma_{k,\perp}^{\prime})^2$ and
$(A^{\prime})^TA^{\prime}V_{\perp}=\mathbf{0}$
as well as $V_k^TV_k^{\perp}=\mathbf{0}, V_{\perp}^TV_k^{\perp}=\mathbf{0}$
and $V_{\perp}^TV_k=\mathbf{0}$, from \eqref{decompqk} we obtain
\begin{align}
&\tilde{q}_k^T(A^{\prime})^TA^{\prime}\tilde{q}_k\notag\\
&=\left(V_k^{\perp}(V_k^{\perp})^T\tilde{q}_k+
V_{\perp}V_{\perp}^T
\tilde{q}_k+V_kV_k^T
\tilde{q}_k\right)^T(A^{\prime})^T A^{\prime}
\left(V_k^{\perp}(V_k^{\perp})^T\tilde{q}_k+V_{\perp}V_{\perp}^T
\tilde{q}_k+
V_kV_k^T\tilde{q}_k\right) \notag\\
&=\left(\tilde{q}_k^TV_k^{\perp}(V_k^{\perp})^T+\tilde{q}_k^TV_{\perp}V_{\perp}^T
+\tilde{q}_k^TV_kV_k^T\right)
\left(V_k^{\perp}(\Sigma_{k,\perp}^{\prime})^2(V_k^{\perp})^T\tilde{q}_k+
V_k(\Sigma_k^{\prime})^2V_k^T
\tilde{q}_k\right) \notag\\
&=\tilde{q}_k^TV_k^{\perp}(\Sigma_{k,\perp}^{\prime})^2(V_k^{\perp})^T\tilde{q}_k
+\tilde{q}_k^TV_k(\Sigma_k^{\prime})^2V_k^T\tilde{q}_k. \label{expansion}
\end{align}
Observe that it is impossible for $(V_k^{\perp})^T\tilde{q}_k$ and
$V_k^T\tilde{q}_k$ to be the eigenvectors of $(\Sigma_{k,\perp}^{\prime})^2$
and $(\Sigma_k^{\prime})^2$ associated with their respective smallest eigenvalues
$\sigma_s^2$ and $\sigma_k^2$ simultaneously, which are
the $(s-k)$-th canonical vector $e_{s-k}^{(s-k)}$ of $\mathbb{R}^{s-k}$ and
the $k$-th canonical vector $e_k^{(k)}$ of $\mathbb{R}^{k}$, respectively;
otherwise, we have $\tilde{q}_k=v_s$ and
$\tilde{q}_k=v_k$ simultaneously, which are impossible as $k<s$. Therefore,
from \eqref{expansion}, \eqref{qktilde} and \eqref{angle2},
we obtain the strict inequality
\begin{align*}
\tilde{q}_k^T(A^{\prime})^TA^{\prime}\tilde{q}_k&> \|(V_k^{\perp})^T\tilde{q}_k\|^2
\sigma_s^2+\|V_k^T\tilde{q}_k\|^2\sigma_k^2
=(1-\varepsilon_k^2)\sigma_s^2+\varepsilon_k^2 \sigma_k^2,
\end{align*}
from which it follows that the lower bound of \eqref{rqi} holds.

Similarly, from \eqref{expansion}, \eqref{qktilde} and \eqref{angle2}
we obtain the upper bound of \eqref{rqi}:
$$
\tilde{q}_k^T(A^{\prime})^TA^{\prime}\tilde{q}_k <\|(V_k^{\perp})^T\tilde{q}_k\|^2
\|(\Sigma_k^{\perp})^2\|+\|V_k^T\tilde{q}_k\|^2\|(\Sigma_k^{\prime})^2\|
=(1-\varepsilon_k^2)\sigma_{k+1}^2+\varepsilon_k^2 \sigma_1^2.
$$

From the lower bound of \eqref{rqi}, we see that if
$\varepsilon_k$ satisfies $\varepsilon_k^2 \sigma_k^2\geq \sigma_{k+1}^2$,
i.e., $\varepsilon_k\geq \frac{\sigma_{k+1}}{\sigma_k}$,
then $\sqrt{\tilde{q}_k^T(A^{\prime})^TA^{\prime}\tilde{q}_k}>\sigma_{k+1}$, i.e.,
\eqref{est1} holds.

Recall from Section~\ref{lsqr} that Algorithm~\ref{alg:lb} generates the same
results when applied to $A$ and $A^{\prime}$. Therefore,
from \eqref{newbk}, we obtain
$B_k^TB_k=Q_k^T(A^{\prime})^TA^{\prime}Q_k$.
Note that $(\theta_k^{(k)})^2$ is the smallest eigenvalue
of the symmetric positive definite matrix $B_k^TB_k$.
Therefore, we have
\begin{equation}\label{rqi2}
(\theta_k^{(k)})^2=\min_{\|c\|=1} c^TQ_k^T(A^{\prime})^TA^{\prime}Q_kc=
\min_{q\in \mathcal{V}_k^R,\ \|q\|=1} q^T(A^{\prime})^TA^{\prime}q
=\widehat{q}_k^T(A^{\prime})^TA^{\prime}\widehat{q}_k.
\end{equation}
Therefore, for $\tilde{q}_k$, we have
$$
\theta_k^{(k)}\leq \sqrt{\tilde{q}_k^T(A^{\prime})^TA^{\prime}\tilde{q}_k},
$$
from which it follows from \eqref{rqi} that
$(\theta_k^{(k)})^2<(1-\varepsilon_k^2)\sigma_{k+1}^2+\varepsilon_k^2 \sigma_1^2$.
As a result, for any $\delta>0$, we can choose $\varepsilon_k\geq 0$ such that
$$
(\theta_k^{(k)})^2<(1-\varepsilon_k^2)\sigma_{k+1}^2+\varepsilon_k^2 \sigma_1^2
\leq (1+\delta)\sigma_{k+1}^2,
$$
i.e., \eqref{thetasigma} holds,
solving which for $\varepsilon_k^2$ gives $\varepsilon_k^2\leq\frac{\delta}
{(\frac{\sigma_1}{\sigma_{k+1}})^2-1}$.
\qquad\endproof

\begin{remark}
The author in \cite{jia18b} has given a detailed
analysis on $\|\sin\Theta(\mathcal{V}_k,\mathcal{V}_k^R)\|$ for the three
kinds of ill-posed problems. It turns out that
$\|\sin\Theta(\mathcal{V}_k,\mathcal{V}_k^R)\|$ cannot be
close to one for severely or moderately ill-posed problems with suitable
$\rho>1$ or $\alpha>1$ and $k\leq k_0$, but it generally approaches one for
mildly ill-posed problems or moderately ill-posed problems with
$\alpha>1$ not enough when $k$ is small.
\end{remark}

\begin{remark}
It has been shown in \cite{jia18b} that for severely and moderately
ill-posed problems with suitable
$\rho>1$ and $\alpha>1$, we may have $\theta_k^{(k)}>\sigma_{k+1}$ for
$k\leq k^*$,
and for mildly ill-posed problems and moderately ill-posed problems
with $\alpha>1$ not enough we have $\theta_k^{(k)}<\sigma_{k+1}$
for some $k\leq k^*$.
\end{remark}

An intrinsic disadvantage of Theorem~\ref{initial} is that
it does not give any sufficient conditions on $\rho$ and $\alpha$ that ensures
$\theta_k^{(k)}>\sigma_{k+1}$. In the next section, we present accurate
results on Problems (i)--(iv) stated in the beginning of Section~\ref{sine}.

\section{The rank $k$ approximation $P_{k+1}B_kQ_k^T$ to $A^{\prime}$, the Ritz values
$\theta_i^{(k)}$ and the regularization of LSQR}\label{rankapp}

For the rank $k$ approximation $P_{k+1}B_kQ_k^T$ to $A^{\prime}$
in LSQR, we define
\begin{equation}\label{gammak}
\gamma_k^{\prime}=\|A^{\prime}-P_{k+1}B_kQ_k^T\|,\ k=1,2,\ldots,s-1,
\end{equation}
which measures the accuracy of the rank $k$ approximation $P_{k+1}B_kQ_k^T$
to $A^{\prime}$. The rank
$k$ matrix $P_{k+1}B_kQ_k^T$ is called a near best rank $k$ approximation
to $A^{\prime}$ if it satisfies
\begin{equation}\label{near}
\sigma_{k+1}\leq \gamma_k^{\prime}<\frac{\sigma_k+\sigma_{k+1}}{2},
\end{equation}
that is, $\gamma_k^{\prime}$ lies between $\sigma_{k+1}$ and $\sigma_k$
and is closer to $\sigma_{k+1}$. This definition
has been introduced in \cite{jia18b} and
shown to be irreplaceable in the context of linear discrete ill-posed problems
when considering the approximation behavior of the Ritz values
$\theta_i^{(k)}$ and the corresponding counterparts involved
in the Krylov solvers CGME and LSMR \cite{jia18d}

With the replacement of the index $n$ by $s$,
the following results in \cite[Theorem 3.2]{jia18b} carry over to
the multiple singular value case.

\begin{theorem}\label{main1}
Assume that the discrete Picard condition \eqref{newpicard2} is
satisfied, and let $|L_{k_1}^{(k)}(0)|$ be defined by \eqref{lk}.
Then for $k=1,2,\ldots,s-1$ we have
\begin{equation}\label{final}
  \sigma_{k+1}\leq \gamma_k^{\prime}\leq \sqrt{1+\eta_k^2}\sigma_{k+1}
\end{equation}
with
\begin{equation} \label{const1}
\eta_k\leq
\xi_k \frac{\max_{k+1\leq j\leq s}|u_i^Tb|}{\min_{1\leq i\leq k}|u_i^Tb|}
\left(1+\mathcal{O}(\rho^{-2})\right), \ k=1,2,\ldots,k_0
\end{equation}
for severely ill-posed problems  with $\sigma_i=\mathcal{O}(\rho^{-i})$,
$i=1,2,\ldots,s$ and
\begin{equation}\label{const2}
\eta_k\leq \left\{\begin{array}{ll}
\xi_1\frac{\sigma_1}{\sigma_2}\frac{\max_{2\leq i\leq n}| u_i^Tb|}{| u_1^Tb|}
\sqrt{\frac{1}{2\alpha-1}} & \mbox{ for } k=1, \\
\xi_k\frac{\sigma_k}{\sigma_{k+1}}\frac{\max_{k+1\leq i\leq n}|u_i^Tb|}
{\min_{1\leq i\leq k}| u_i^T b|}
\sqrt{\frac{k^2}{4\alpha^2-1}+\frac{k}{2\alpha-1}}|L_{k_1}^{(k)}(0)|
& \mbox{ for } 1< k\leq k_0
\end{array}
\right.
\end{equation}
for moderately or mildly ill-posed problems with $\sigma_j=\zeta j^{-\alpha},\
j=1,2,\ldots,s$, where
$\xi_k=\sqrt{\left(\frac{\|\Delta_k\|}{1+\|\Delta_k\|^2}\right)^2+1}$ for
$\|\Delta_k\|<1$ and $\xi_k\leq\frac{\sqrt{5}}{2}$ for $\|\Delta_k\|\geq 1$.
\end{theorem}

Based on Theorem~\ref{main1}, for the ill-posed problems with the
singular value models $\sigma_k=\zeta\rho^{-k}$ and $\sigma_k=\zeta k^{-\alpha}$,
the following two theorems establish the sufficient conditions on $\rho$ and
$\alpha$ that guarantee that $P_{k+1}B_kQ_k^T$ is a near best rank $k$
approximation to $A^{\prime}$ and the Ritz values $\theta_i^{(k)}$ approximate
the $k$ large singular values $\sigma_i$ of $A^{\prime}$ in natural
order for $k=1,2,\ldots,k^*$, whose proofs are
the same as those of Theorem 3.3 and Theorem 4.1 in \cite{jia18b}.

\begin{theorem}\label{nearapprox}
For a given \eqref{eq1}, assume that the discrete Picard condition
\eqref{newpicard2} is satisfied. Then, in the sense of \eqref{near},
$P_{k+1}B_kQ_k^T$ is a near best rank $k$ approximation to $A^{\prime}$
for $k=1,2,\ldots,k^*$ if
\begin{equation}\label{condition}
\sqrt{1+\eta_k^2}<\frac{1}{2}\frac{\sigma_k}{\sigma_{k+1}}+\frac{1}{2}.
\end{equation}
Furthermore, $P_{k+1}B_kQ_k^T$ is a near best rank $k$ approximation to $A^{\prime}$
if $\rho>\frac{1+\sqrt{5}}{2}$ for the severely ill-posed problems
with $\sigma_k=\zeta\rho^{-k}$ or $\alpha$ satisfies
\begin{equation}\label{condition1}
2\sqrt{1+\eta_k^2}-1<\left(\frac{k+1}{k}\right)^{\alpha},\
k=1,2,\ldots,k^*
\end{equation}
for the moderately and mildly ill-posed problems with $\sigma_k=\zeta k^{-\alpha}$,
respectively.
\end{theorem}

The author in \cite{jia18a} has given a detailed analysis on this theorem;
see Remarks 3.8-3.9 there. The conclusions are that, for
severely and moderately ill-posed problems with suitable $\rho>1$ and
$\alpha>1$, Algorithm~\ref{alg:lb} always generates near best
rank $k$ approximations for $k=1,2,\ldots,k^*$ but
$P_{k+1}B_kQ_k^T$ may not be a near best rank $k$ approximation for some
$k\leq k^*$ for moderately ill-posed problems
with $\alpha>1$ not enough and mildly ill-posed problems.

\begin{theorem}\label{ritzvalue}
Assume that \eqref{eq1} is severely ill-posed with
$\sigma_i=\zeta\rho^{-i}$ and $\rho>1$ or moderately ill-posed
with $\sigma_i=\zeta i^{-\alpha}$ and $\alpha>1$, $i=1,2,\ldots,s$,
and the discrete Picard condition \eqref{newpicard2} is
satisfied. Let the Ritz values $\theta_i^{(k)}$ be labeled
as $\theta_1^{(k)}>\theta_2^{(k)}>\cdots>\theta_{k}^{(k)}$.
Then
\begin{align}
\sigma_i-\theta_i^{(k)} &\leq \gamma_k^{\prime}\leq \sqrt{1+\eta_k^2}\sigma_{k+1},\
i=1,2,\ldots,k.\label{error}
\end{align}
For $k=1,2,\ldots,k^*$, if $\rho\geq 1+\frac{\sqrt{6}}{2}$ or $\alpha>1$ satisfies
\begin{equation}\label{condm}
1+\sqrt{1+\eta_{k}^2}<\left(\frac{k+1}{k}\right)^{\alpha},
\end{equation}
then the $k$ Ritz values $\theta_i^{(k)}$ strictly interlace
the first large $k+1$ singular values of $A^{\prime}$ and approximate
the first $k$ large ones in natural order:
\begin{align}
\sigma_{i+1}&<\theta_i^{(k)}<\sigma_i,\,i=1,2,\ldots,k.
\label{error2}
\end{align}
\end{theorem}

\begin{remark}
From Theorem~\ref{nearapprox} and Theorem~\ref{ritzvalue},
it is known that the near best rank $k$ approximation
$P_{k+1}B_kQ_k^T$ to $A$ essentially means that the singular values
$\theta_i^{(k)}$ of $B_k$ approximate the
large singular values $\sigma_i$ of $A^{\prime}$ in natural order for suitable
$\rho>1$ and $\alpha>1$. On the other hand, for a {\em given} problem
with $\alpha>1$, the smaller $k$ is, the more likely
the $k$ Ritz values $\theta_i^{(k)}$ to
approximate the large singular values of $A^{\prime}$ in natural order. Hence
the $k$ Ritz values $\theta_i^{(k)}$ may not
approximate the large singular values of $A^{\prime}$ in natural order
at some $k\leq k^*$ for $\alpha>1$ not enough; in
this case, by Theorem~\ref{propaprime} we must have $k^*<k_0$.
\end{remark}

\begin{remark}\label{extract}
Theorem~\ref{nearapprox} and Theorem~\ref{ritzvalue} show that LSQR has the full
regularization for these two kinds of ill-posed
problems with suitable $\rho>1$ and $\alpha>1$.
\end{remark}

\begin{remark}
For mildly ill-posed
problems, we observe that the sufficient condition \eqref{condm} for \eqref{error2}
is never met because
$
\left(\frac{k+1}{k}\right)^{\alpha}\leq 2
$
for any $k\geq 1$ and $\frac{1}{2}< \alpha\leq 1$. This indicates that
the $k$ Ritz values $\theta_i^{(k)}$ may not approximate
the large singular values $\sigma_i$ of $A^{\prime}$ in natural order
soon as $k$ increases.
\end{remark}

\section{Monotonicity of $\gamma_k^{\prime}$ and
decay rates of the entries $\alpha_k$ and $\beta_{k+1}$}
\label{alphabeta}

In this section, we extend the results on the monotonicity of $\gamma_k^{\prime}$
and the decay rates of $\alpha_k$ and $\beta_{k+1}$ in \cite{jia18b} to the
multiple singular value case by making some changes in the proofs.

\begin{theorem}\label{main2}
With the notation defined previously, the following results hold:
\begin{eqnarray}
  \sqrt{\alpha_{k+1}^2+\beta_{k+2}^2}
  &<&\gamma_k^{\prime}\leq \sqrt{1+\eta_k^2}\sigma_{k+1},
  \ k=1,2,\ldots,s-1,\label{alpha}\\
 \alpha_{k+1}\beta_{k+2}&\leq &
\frac{(\gamma_k^{\prime})^2}{2}\leq
   \frac{(1+\eta_k^2)\sigma_{k+1}^2}{2}, \ k=1,2,\ldots,s-1,
  \label{prod2}\\
\gamma_{k+1}^{\prime}&<&\gamma_k^{\prime},\  \ k=1,2,\ldots,s-2.
\label{gammamono}
\end{eqnarray}
\end{theorem}

{\em Proof}.
In the multiple singular value case, as we have shown
in Section~\ref{lsqr}, Algorithm~\ref{alg:lb}
can only be run to step $s$ without breakdown. From \eqref{lbaprime}
we augment $P_{s+1}$ and $Q_s$ to the $m\times m$ and $n\times n$
orthogonal matrices $P=(P_{s+1},\hat{P})$ and $Q=(Q_s,\hat{Q})$,
respectively. Then from \eqref{lbaprime} we obtain
$$
P^TAQ=P^TA^{\prime}Q=\left(\begin{array}{cc}
B_s&\mathbf{0}\\
\mathbf{0}& \mathbf{0}
\end{array}
\right),
$$
from which it follows that
$$
\gamma_k^{\prime}=\|A^{\prime}-P_{k+1}B_kQ_k^T\|
=\|P^T\left(A^{\prime}-P_{k+1}B_kQ_k^T\right)Q\|=\|G_k^{\prime}\|,
$$
where $G_k^{\prime}$ is the right bottom $(s-k+1)\times (s-k)$ matrix of
$B_s$. Then the rest proof is exactly the same as that of the
Theorem 5.1 in \cite{jia18b}.
\qquad\endproof

\begin{remark}
The strict decreasing property \eqref{gammamono} of $\gamma_k^{\prime}$ and
the lower bounds on $\gamma_k^{\prime}$ in \eqref{alpha}--\eqref{prod2} hold
unconditionally for a general $A$,
independent of the degree of ill-posedness.
\end{remark}

\begin{remark}
It is impractical to compute $\gamma_k^{\prime}$ for $A$ large.
However, \eqref{alpha} and \eqref{prod2} indicates that the sum
$\alpha_{k+1}+\beta_{k+2}\leq \sqrt{2}\gamma_k^{\prime}$ and
decays as fast as $\gamma_k^{\prime}$. Therefore,
we can reliably judge the decay rates of $\gamma_k^{\prime}$
during computation with little extra cost.
\end{remark}

\begin{remark}
For the severely and moderately ill-posed problems with suitable $\rho>1$ and
$\alpha>1$,
\eqref{alpha} and \eqref{prod2} show that $\alpha_{k+1}+\beta_{k+2}$
decays as fast as $\sigma_{k+1}$ for $k\leq k_0$. For mildly ill-posed problems,
since $\eta_k$ are generally bigger than one considerably as $k$ increases,
$\alpha_{k+1}+\beta_{k+2}$ cannot generally decay as fast as
$\sigma_{k+1}$.
\end{remark}

\section{Best, near best and general rank $k$ approximations to $A$ and their
implications on LSQR and some others}
\label{morerank}

In this section, we discuss some important issues on
best, near best and general rank $k$ approximations to $A$
when all the singular values of $A$ are simple, i.e.,
the multiplicities $c_i=1,\ i=1,2,\ldots,s$ and $s=n$.
The issues to be addressed has received no attention
in \cite{jia18a,jia18b} and the literature. If $A$ has at least
one multiple singular value, i.e., some $c_i>1$, we speak of
the corresponding rank $k$ approximations to $A^{\prime}$. Therefore,
without loss of generality, in this section
we assume that $A^{\prime}=A,\Sigma^{\prime}=\Sigma
={\rm diag}(\sigma_1,\sigma_2,\ldots,\sigma_n)$,
$U=(u_1,u_2,\ldots,u_n)$ and $V=(v_1,v_2,\ldots,v_n)$
in \eqref{aprime}, and the compact SVD of $A$ is $A=U\Sigma V^T$.
Denote by
\begin{equation}\label{lsqrapp}
\gamma_k=\|A-P_{k+1}B_kQ_k^T\|,\ k=1,2,\ldots,n-1
\end{equation}
the accuracy of the rank $k$ approximation $P_{k+1}B_kQ_k^T$ to $A$.

We first investigate general best or near best rank $k$ approximations
to $A$ with $\sigma_k=\zeta k^{-\alpha}$ and $\alpha>0$.
We will show that, for each of such rank
$k$ approximations, its smallest
nonzero singular value may be smaller than $\sigma_{k+1}$ for
$\alpha\leq 1$, that is, its nonzero singular values
may not approximate the $k$ large singular values of $A$ in natural
order, but it is guaranteed to be bigger
than $\sigma_{k+1}$ for suitable $\alpha>1$.
The implication is that a {\em general} best or near best rank $k$ approximation
to $A$ may have very small nonzero singular values
and thus may {\em not} be a suitable replacement of $A$ for
$\alpha\leq 1$. We then consider the properties of
the Ritz values $\theta_i^{(k)},\ i=1,2,\ldots,k$ and derive
some interesting properties on them when $P_{k+1}B_kQ_k^T$
is or is not a near best rank $k$ approximation to $A$ for
the ill-posed problems with $0<\alpha \leq 1$, which include
mildly ill-posed ones. Finally, we elaborate how to apply
these properties to several other Krylov solvers that solve
ill-posed problems.

First of all, we mention an intrinsic fact that
best rank $k$ approximations to $A$ with respect to the 2-norm are not
{\em unique}. In fact, besides $A_k=U_k\Sigma_kV_k^T$ with
$U_k=(u_1,u_2,\ldots,u_k)$, $V_k=(v_1,v_2,\ldots,v_k)$ and
$\Sigma_k={\rm diag}(\sigma_1,\sigma_2,\ldots,\sigma_k)$,
there are other infinitely many best rank $k$ approximations
to $A$. This is certainly also true for near best rank $k$
approximations to $A$, as we see below.

Let $C_k$ be a best or near best rank
$k$ approximation to $A$ with $\|A-C_k\|=(1+\epsilon)\sigma_{k+1}$ with any
$\epsilon\geq 0$ satisfying $(1+\epsilon)\sigma_{k+1}<\frac{\sigma_k+
\sigma_{k+1}}{2}$, that is, $(1+\epsilon)\sigma_{k+1}$ is between $\sigma_{k+1}$
and $\sigma_k$ and closer to $\sigma_{k+1}$ (Note: $\epsilon=0$ corresponds to
a best rank $k$ approximation $C_k$.). Then we have
\begin{equation}\label{nearrank}
1+2\epsilon<\frac{\sigma_k}{\sigma_{k+1}}.
\end{equation}
It is remarkable to note that $C_k$ is not unique for any
given $\epsilon\geq 0$ satisfying \eqref{nearrank}. For example, among
others, it is easy to verify that all the matrices
\begin{eqnarray}
C_k&=&A_k(\theta,j)\nonumber\\
&=&A_k-\sigma_{k+1}U_k
{\rm diag}(\theta(1+\epsilon),\ldots,\theta(1+\epsilon),
\underbrace{(1+\epsilon)}_j,\theta(1+\epsilon),\ldots,\theta(1+\epsilon))
V_k^T \label{ck}
\end{eqnarray}
with any $\theta \in [0,1]$ and $1\leq j\leq k-1$ form a family of
best or near best rank $k$ approximations to $A$ that satisfy
$\|A-C_k\|=(1+\epsilon)\sigma_{k+1}$. Meanwhile, it is easily
seen that the smallest nonzero singular value of $C_k$ is
$\sigma_k-\theta(1+\epsilon)\sigma_{k+1}$.

\begin{theorem}\label{nearappro}
Let $C_k$ be the best or near best approximations to $A$ defined as
\eqref{ck}. Then
for $0<\alpha\leq 1$, if $\theta$ is
sufficiently close to one and $k>1$, it holds that
\begin{equation}\label{misinter}
\sigma_k-\theta(1+\epsilon)\sigma_{k+1}<\sigma_{k+1},
\end{equation}
that is, the smallest nonzero singular value
$\sigma_k-\theta(1+\epsilon)\sigma_{k+1}$ of $C_k$ does not lie between
$\sigma_{k+1}$ and $\sigma_k$; if $\theta$ is sufficiently close to zero,
then $\sigma_k-\theta(1+\epsilon)\sigma_{k+1}$ lies in
$\sigma_{k+1}$ and $\sigma_k$:
\begin{equation}\label{inter}
\sigma_{k+1}<\sigma_k-\theta(1+\epsilon)\sigma_{k+1}<\sigma_k.
\end{equation}
If $\alpha>1$ sufficiently, \eqref{inter} holds for $k>1$
and any $\theta\in [0,1]$.
\end{theorem}

{\em Proof}.
Since $\sigma_k=\zeta k^{-\alpha}$ and
$
\left(\frac{k+1}{k}\right)^{\alpha}<2
$
for any $k>1$ and $0<\alpha\leq 1$,
we obtain
\begin{equation}\label{mildmoderate2}
\sigma_k-\theta(1+\epsilon)\sigma_{k+1}=\sigma_{k+1}
\left(\left(\frac{k+1}{k}\right)^{\alpha}-\theta(1+\epsilon)\right)<\sigma_{k+1}
\end{equation}
for $\theta$ sufficiently close to one.
This shows that $\sigma_k-\theta(1+\epsilon)\sigma_{k+1}$ does not lie between
$\sigma_{k+1}$ and $\sigma_k$ and does not interlace them for $k>1$.
In this case, for a given $\alpha\leq 1$, the bigger $k$ is, the
smaller $\left(\frac{k+1}{k}\right)^{\alpha}-\theta(1+\epsilon)$ is, and the
further is $\sigma_k-\theta(1+\epsilon)\sigma_{k+1}$ away from $\sigma_{k+1}$.
On the other hand, for $\theta$ sufficiently small we
have
\begin{equation}\label{mildmoderate}
\left(\frac{k+1}{k}\right)^{\alpha}-\theta(1+\epsilon)>1,
\end{equation}
that is, $\sigma_k-\theta(1+\epsilon)\sigma_{k+1}$ interlaces
$\sigma_{k+1}$ and $\sigma_k$ for $\theta$
sufficiently small.

For $A$ with $\sigma_k=\zeta k^{-\alpha}$ and $\alpha>1$ and
$k=1,2,\ldots,n-1$, the requirement \eqref{mildmoderate} is met
for {\em any} $\theta\in [0,1]$ and suitable $\alpha>1$,
leading to $\sigma_k-\theta(1+\epsilon)\sigma_{k+1}>\sigma_{k+1}$. This means
that the smallest nonzero singular value
$\sigma_k-\theta(1+\epsilon)\sigma_{k+1}$ of the best or near best rank
approximation $C_k$ lies between $\sigma_{k+1}$ and $\sigma_k$ for
suitable $\alpha>1$.
\qquad\endproof

We should be aware that \eqref{misinter} is established by assuming
the {\em worst} case that, over all best or
near best rank $k$ approximations $C_k$ of form \eqref{ck},
the minimum $\sigma_k-(1+\epsilon)\sigma_{k+1}$
of the smallest nonzero singular values of the $C_k$ is almost or exactly
taken, i.e., $\theta\approx 1$ or $\theta=1$.
We now prove that the minimum is indeed
$\sigma_k-(1+\epsilon)\sigma_{k+1}$ over the set of {\em all} the
near best rank $k$ approximations $C_k$, including the ones of form \eqref{ck}.
Let $\sigma_k(C_k)$ be the smallest nonzero singular value of $C_k$.
Then from $\|A-C_k\|=(1+\epsilon)\sigma_{k+1}$,
by the standard perturbation theory we have
$$
|\sigma_k-\sigma_k(C_k)|\leq (1+\epsilon)\sigma_{k+1}.
$$
Clearly, the minimum of all the $\sigma_k(C_k)$ is attained if and
only if the above equality holds and the left-hand side is positive,
which means that this minimum is exactly $\sigma_k-(1+\epsilon)\sigma_{k+1}$.
In contrast, \eqref{inter} holds essentially in the {\em best} case that
the maximum of the smallest singular values
$\sigma_k-\theta(1+\epsilon)\sigma_{k+1}$ of $C_k$ defined by
\eqref{ck} is almost or exactly taken, i.e., $\theta\approx 0$ or
$\theta=0$.

As far as LSQR is concerned, notice from \cite{jia18b} that
condition \eqref{condm} for the
interlacing property \eqref{error2} is derived by assuming the worst case that
$\sigma_k-\theta_k^{(k)}=\gamma_k\leq\sqrt{1+\eta_k^2}\sigma_{k+1}$, that is,
$\theta_k^{(k)}$ is supposed to be the smallest possible nonzero one among
all the $\sigma_k(C_k)$, where $C_k$ belongs to the set of near best
$k$ approximations that satisfy $\|A-C_k\|=\gamma_k\leq\sqrt{1+\eta_k^2}
\sigma_{k+1}$. Even so,
a near best rank $k$ approximation $P_{k+1}B_kQ_k^T$ can guarantee
the approximations of $\theta_i^{(k)},\ i=1,2,\ldots,k$
to the large singular values $\sigma_i$ in natural order for suitable $\alpha>1$.
This is in accordance
with Theorems~\ref{nearapprox}--\ref{ritzvalue} though the sizes of $\alpha>1$
for them are different.

For mildly ill-posed problems, Theorem~\ref{nearappro} and
the above analysis indicate that the $k$ Ritz values $\theta_i^{(k)}$
may or may not approximate the large singular values $\sigma_i$
of $A$ in natural order {\em even} when $P_{k+1}B_kQ_k^T$ is a near best
rank $k$ approximation to $A$.

Unfortunately, as we have elaborated and numerically confirmed in \cite{jia18b},
for mildly ill-posed problems $P_{k+1}B_kQ_k^T$ may be a near best rank
$k$ approximation to $A$ only for $k$ very small and it is
not any more soon as $k$ increases. The following results are
more general and cover all ill-posed problems with $0<\alpha\leq 1$,
which include mildly ill-posed ones. We will seek the maximum possible
number of the Ritz values smaller than $\sigma_{k+1}$ and
get insight into how small they can be.

\begin{theorem}\label{disorder}
For $\sigma_i=\zeta i^{-\alpha},\,i=1,2,\ldots,n$ with $0<\alpha\leq 1$,
suppose $\gamma_k\in [\sigma_{j+1},\sigma_j]$ for some $j\leq k$.
Then
\begin{enumerate}
\item if $\gamma_k$ is closer to $\sigma_{j+1}$ than to $\sigma_j$
and $j(j+1)\geq k$,  there are at most $k-j+1$ Ritz values
$\theta_j^{(k)},\theta_{j+1}^{(k)},\ldots,\theta_k^{(k)}$ smaller
than $\sigma_{k+1}$;

\item for $j\geq 2$,
if $\gamma_k$ is closer to $\sigma_j$ than to $\sigma_{j+1}$
and $j(j+1)>k$, there are at most $k-j+2$ Ritz values
$\theta_{j-1}^{(k)},\theta_j^{(k)},\ldots,\theta_k^{(k)}$ smaller
than $\sigma_{k+1}$;

\item for $j=1$, if $\gamma_k$ is closer to $\sigma_1$ than to $\sigma_2$,
all the Ritz values $\theta_i^{(k)},\ i=1,2,\ldots,k$ are
possibly smaller than $\sigma_{k+1}$.
\end{enumerate}
\end{theorem}
{\em Proof}.
Since $\sigma_j=\zeta j^{-\alpha}$, we have
\begin{align}
\frac{\alpha}{j}\sigma_{j+1}&=\frac{\alpha}{\zeta j^{1-\alpha}}\zeta
j^{-\alpha}\sigma_{j+1}
=\frac{\alpha}{\zeta j^{1-\alpha}}\zeta^2 j^{-\alpha} (j+1)^{-\alpha}
\nonumber\\
&=
\frac{\alpha}{\zeta j^{1-\alpha}}\zeta \sigma_{j(j+1)}=
\frac{\alpha}{j^{1-\alpha}}\sigma_{j(j+1)}. \label{step}
\end{align}
Define the function $f(x)=(1+x)^{\alpha}$ for
$x>0$ and
$0<\alpha\leq 1$. Then its $\ell$th derivative
$$
f^{(\ell)}(x)=\alpha(\alpha-1)\cdots (\alpha-\ell+1)(1+x)^{(\alpha-1)(\alpha-2)
\cdots (\alpha-\ell)}
$$
is always negative for $\ell\geq 2$. Therefore, taking the first {\em three}
terms of the Taylor expansion of $f(\frac{1}{j})=(1+\frac{1}{j})^{\alpha}$
and exploiting \eqref{step}, we obtain
\begin{eqnarray}
\sigma_j-\sigma_{j+1}&=&
\sigma_{j+1}\left(\left(\frac{j+1}{j}\right)^{\alpha}-1\right)
=\sigma_{j+1}\left(\left(1+\frac{1}{j}\right)^{\alpha}-1\right)\nonumber\\
&\geq& \sigma_{j+1}\left(1+\frac{\alpha}{j}-\frac{(1-\alpha)\alpha}
{2j^2}-1\right) \nonumber\\
&=&\frac{\alpha}{j}\sigma_{j+1} \left(1-\frac{(1-\alpha)}{2j}\right)\nonumber\\
&=&\frac{\alpha}{j^{1-\alpha}}\sigma_{j(j+1)}\left(1-\frac{(1-\alpha)}{2j}\right).
\label{step2}
\end{eqnarray}

In the case that $\gamma_k$ lies in $[\sigma_{j+1},\sigma_j]$ and
is closer to $\sigma_{j+1}$ than to $\sigma_j$  for $1\leq j\leq k$,
by taking $i=j$ in
\eqref{error}, we get $\sigma_j-\theta_j^{(k)}\leq\gamma_k$.
Therefore, by the assumption on $\gamma_k$,
from \eqref{step2} we obtain the lower bounds for $\theta_j^{(k)}$:
\begin{equation}\label{sharp}
\theta_j^{(k)}\geq\sigma_j-\gamma_k\geq \frac{\sigma_j-\sigma_{j+1}}{2}\geq
\frac{\alpha}{2j^{1-\alpha}}\sigma_{j(j+1)}\left(1-\frac{(1-\alpha)}{2j}\right).
\end{equation}
Since each lower bound in \eqref{sharp} cannot be improved.
$\theta_j^{(k)}$ is always likely to (approximately) attain the second lower bound.
Suppose it is the case, and note that $\frac{\alpha}{j^{1-\alpha}}<1$
decreases with $j$ for a given $\alpha \in (0,1)$ and equals one
for $\alpha=1$. We then obtain
\begin{eqnarray*}
\theta_j^{(k)}&\approx &
\frac{\alpha}{2j^{1-\alpha}}\sigma_{j(j+1)}\left(1-\frac{(1-\alpha)}{2j}\right)\\
&\leq& \frac{\alpha}{2j^{1-\alpha}}\sigma_{j(j+1)}\leq
\frac{1}{2}\sigma_{j(j+1)}<\sigma_{j(j+1)},
\end{eqnarray*}
which is smaller than $\sigma_{k+1}$ when $j(j+1)\geq k$.

Moreover, whenever $\theta_j^{(k)}<\sigma_{k+1}$,
by the labeling rule, there are $k-j+1$ Ritz values
$\theta_j^{(k)},\theta_{j+1}^{(k)},\ldots,\theta_k^{(k)}$ smaller
than $\sigma_{k+1}$.

We next analyze the case that $\gamma_k$ lies in $[\sigma_{j+1},\sigma_j]$ and
is closer to $\sigma_j$ than
to $\sigma_{j+1}$ for $j\geq 2$. Taking $i=j-1$ in \eqref{error}, we have
$
\sigma_{j-1}-\theta_{j-1}^{(k)}\leq \gamma_k,
$
which shows that
$$
\theta_{j-1}^{(k)}\geq \sigma_{j-1}-\gamma_k.
$$
For $\sigma_j=\zeta j^{-\alpha}$ with $\alpha>0$, it is easy to justify that
$$
\sigma_{j-1}-\sigma_j>\sigma_j-\sigma_{j+1}.
$$
Therefore, from \eqref{step2} we obtain
\begin{align*}
\theta_{j-1}^{(k)}&\geq \sigma_{j-1}-\gamma_k\geq \sigma_{j-1}-\sigma_j> \sigma_j-\sigma_{j+1}\\
&\geq \frac{\alpha}{j^{1-\alpha}}\sigma_{j(j+1)}\left(1-\frac{(1-\alpha)}{2j}\right).
\end{align*}

Analogously, when the above lower bound is (approximately) attainable, for $j\geq 2$
we have
\begin{eqnarray*}
\theta_{j-1}^{(k)}&\approx &
\frac{\alpha}{j^{1-\alpha}}\sigma_{j(j+1)}\left(1-\frac{(1-\alpha)}{2j}\right)\\
&\leq& \frac{\alpha}{j^{1-\alpha}}\sigma_{j(j+1)}\leq \sigma_{j(j+1)},
\end{eqnarray*}
which is smaller than $\sigma_{k+1}$ when $j(j+1)>k$.

Whenever $\theta_{j-1}^{(k)}<\sigma_{k+1}$,
by the labeling rule, there are $k-j+2$ Ritz values
$\theta_{j-1}^{(k)},\theta_j^{(k)},\ldots,\theta_k^{(k)}$ smaller
than $\sigma_{k+1}$.

The above analysis does not include the case that $j=1$ and $\gamma_k$ is
closer to $\sigma_1$ than to $\sigma_2$, which needs a special treatment.
For this case, taking $i=1$ in \eqref{error} yields
$
\sigma_1-\theta_1^{(k)}\leq\gamma_k,
$
from which we obtain
\begin{equation}\label{j1}
\theta_1^{(k)}\geq \sigma_1-\gamma_k.
\end{equation}
Under the assumption, such lower bound can be arbitrarily small, which means
that $\theta_1^{(k)}$ can be arbitrarily small so that it can be smaller than
$\sigma_{k+1}$. The closer $\gamma_k$ is to $\sigma_1$, the more likely
$\theta_1^{(k)}<\sigma_{k+1}$. As a consequence, it is possible
that all the Ritz values $\theta_i^{(k)},\
i=1,2,\ldots,k$ are smaller than $\sigma_{k+1}$.
\qquad\endproof

This theorem indicates that, whenever
$j(j+1)>k$ substantially, $\theta_j^{(k)} \mbox{ or $\theta_{j-1}^{(k)}$}
\approx \sigma_{j(j+1)}$ can be considerably smaller
than $\sigma_{k+1}$. Furthermore, we see from \eqref{step} and the proof
that for $j\geq 2$, the smaller $\alpha<1$, the smaller
$\frac{\alpha}{j^{1-\alpha}}\sigma_{j(j+1)}$ is than $\sigma_{j(j+1)}$,
and, for a fixed $\alpha<1$,
the bigger $j$, the smaller $\frac{\alpha}{j^{1-\alpha}}\sigma_{j(j+1)}$
than $\sigma_{j(j+1)}$. Consequently,
the smaller $\alpha<1$, the more likely is
$\theta_j^{(k)}$ or $\theta_{j-1}^{(k)}$ smaller than $\sigma_{k+1}$ when
$j(j+1)>k$. Even for $j=k$, this theorem illustrates that
a near best approximation $P_{k+1}B_kQ_k^T$ does not guarantee that
$\theta_i^{(k)},\ i=1,2,\ldots, k$ approximate
the large singular values $\sigma_i$ of $A$ in natural order;
in the worst case, we may have $\theta_k^{(k)}\leq \frac{\alpha}{k^{1-\alpha}}
\sigma_{k(k+1)}<\sigma_{k(k+1)}$.

More generally, if $P_{k+1}B_kQ_k^T$ is replaced by
any rank $k$ approximation $C_k$ to $A$ and we still denote
$\gamma_k=\|A-C_k\|$ and by $\theta_i^{(k)}$
the $k$ nonzero singular values of $C_k$ by requiring
that $\theta_i^{(k)}\leq \sigma_i$, $i=1,2,\ldots,k$,
then this theorem still holds.
Therefore, for such a $C_k$ and $0<\alpha\leq 1$,
the nonzero singular values of $C_k$ do not necessarily approximate
the large singular values of $A$ in natural order, and some of them may
be smaller than $\sigma_{k+1}$.

MINRES and MR-II are Krylov solvers for
\eqref{eq1} with $A$ symmetric and have been shown to have
regularizing effects \cite{hanke95,hanke96,hansen10,huangjia17,jensen07,kilmer99},
but MR-II is preferable since
the noisy $b$ is excluded in the underlying subspace \cite{huangjia17,jensen07}.
For $A$ nonsymmetric or multiplication with $A^{T}$ difficult to compute,
GMRES and RRGMRES are candidate methods, and the latter
may be better \cite{jensen07}. We mention that the regularizing
effects of GMRES and RRGMRES are highly problem
dependent, and it appears that they require that the mixing of the left and
right singular vectors of $A$ be weak, that is, $V^TU$ is
close to a diagonal matrix; see,
e.g., \cite{jensen07} and \cite[p.126]{hansen10}.

Similar to LSQR, all the above methods and
CGME \cite{bjorck96,bjorck15,craig,hanke95,hanke01} generate their own rank $k$
approximations to $A$ at iteration $k$.
For each of these rank $k$ approximation to $A$,
still denote by $\theta_i^{(k)},\ i=1,2,\ldots,k$
its nonzero singular values. Since the first inequality
of \eqref{error} holds for all these methods. As a consequence,
Theorem~\ref{disorder} holds for all of them, and its proof carries
over to these methods without any change.

However, we must point out that, for a given \eqref{eq1},
the size and properties of $\gamma_k$ in these methods
differ greatly, so is the approximation
behavior of $\theta_i^{(k)},\ i=1,2,\ldots,k$. For example, $\gamma_k$ in CGME
monotonically decreases strictly with respect to $k$, and is
unconditionally bigger than that in LSQR \cite{jia18d}. Though not explicitly
presented in \cite{huangjia17}, using the analysis approach
in \cite{jia18b,jia18d}, we can easily prove that $\gamma_k$ in MINRES or MR-II
monotonically decreases with respect to $k$. Such monotonic decreasing
property of $\gamma_k$ has turned out to be extremely important
for a Krylov solver in the context of ill-posed problems since the rank
$k$ approximation is becoming increasingly better replacement of
$A$ when $gamma_k$ decreases monotonically.
If the Krylov solver does not have this property,
it may not be a good regularization method
for solving \eqref{eq1}. Unfortunately, for a general nonsymmetric $A$,
it can be justified that the $\gamma_k$ in GMRES and RRGMRES,
mathematically, do not have such monotonic decreasing property,
and can behave very irregularly with $k$ increasing. These mean
that these two methods could not have regularizing effects
and may not find any meaningful regularized solutions.

\section{Numerical experiments}\label{exper}

In this section, we present numerical experiments to
justify Theorems~\ref{nearapprox}--\ref{ritzvalue}, Theorem~\ref{main2}
and Theorems~\ref{nearappro}--\ref{disorder} when $A$ is supposed to have only
simple singular values. The extensive
experiments in \cite{jia18a,jia18b,jia18d} have confirmed the other
results when the singular values $\sigma_i$ of $A$ are all simple. To this end,
we use some random matrices {\sf regutm} and the Matlab
code {\sf regutm.m} from \cite{hansen07} to generate the test problems. The
function {\sf [A,U,V]=regutm(m,n,s)}
constructs a random $m \times n$ matrix $A$ in the normal distribution
such that the eigenvectors $u_k$ and $v_k$
of $AA^T$ and $A^TA$ are oscillating. Precisely, it generates
$$
A = U {\rm diag}(s) V^T,
$$
where the number of sign changes in $u_k$ and $v_k$ is exactly $k-1$, and
the third argument $s$ specifies the singular values $\sigma_k$ of $A$.

In our experiments, we first take
$\sigma_k=\frac{1}{k^{\alpha}},\ k=1,2,\ldots,n$
by taking $m=n=10,000$ and $\alpha=1,0.6$ and $0.3$, respectively,
which meet the assumptions of Theorems~\ref{nearappro}--\ref{disorder}.
We construct the exact solution $x_{true}=ones(n,1)$ and the noise free
right-hand side $b_{true}=Ax_{true}$. Then we generate a Gaussian white
noise vector $e$ with zero mean and the relative noise level
$\varepsilon=\frac{\|e\|}{\|b_{true}\|}=10^{-3}$, and
add it to $b_{true}$ to form the noisy right-hand side $b$.
We use the LSQR algorithm with the
starting vector $p_1=b/\|b\|$ to solve $Ax=b$.
Here we are {\em only} concerned
with the accuracy $\gamma_k$ of the rank $k$ approximations $P_{k+1}B_kQ_k^T$
and the approximation properties of the Ritz values $\theta_i^{(k)},
\ i=1,2,\ldots,k$. To simulate exact arithmetic, we use Algorithm~\ref{alg:lb}
with reorthogonalization to generate $B_k$ and numerically
orthonormal $P_{k+1}$ and $Q_k$.

All the computations are carried out in Matlab R2017b on the
Intel Core i7-4790k with CPU 4.00 GHz processor and 16 GB RAM
with the machine precision
$\epsilon_{\rm mach}= 2.22\times10^{-16}$ under the Miscrosoft
Windows 8 64-bit system.

\begin{figure}
\begin{minipage}{0.48\linewidth}
  \centerline{\includegraphics[width=6.0cm,height=4.5cm]{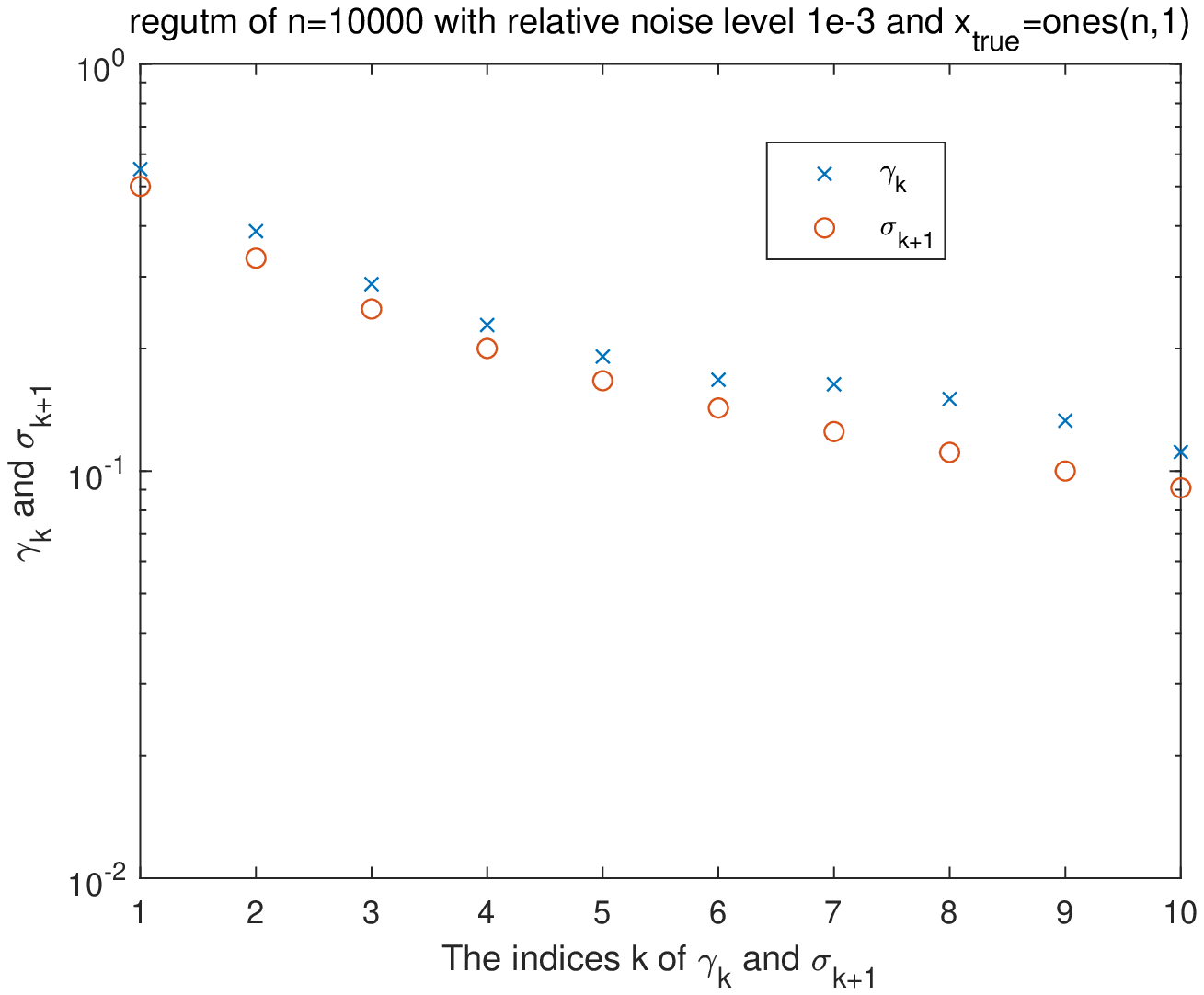}}
  \centerline{(a)}
\end{minipage}
\hfill
\begin{minipage}{0.48\linewidth}
  \centerline{\includegraphics[width=6.0cm,height=4.5cm]{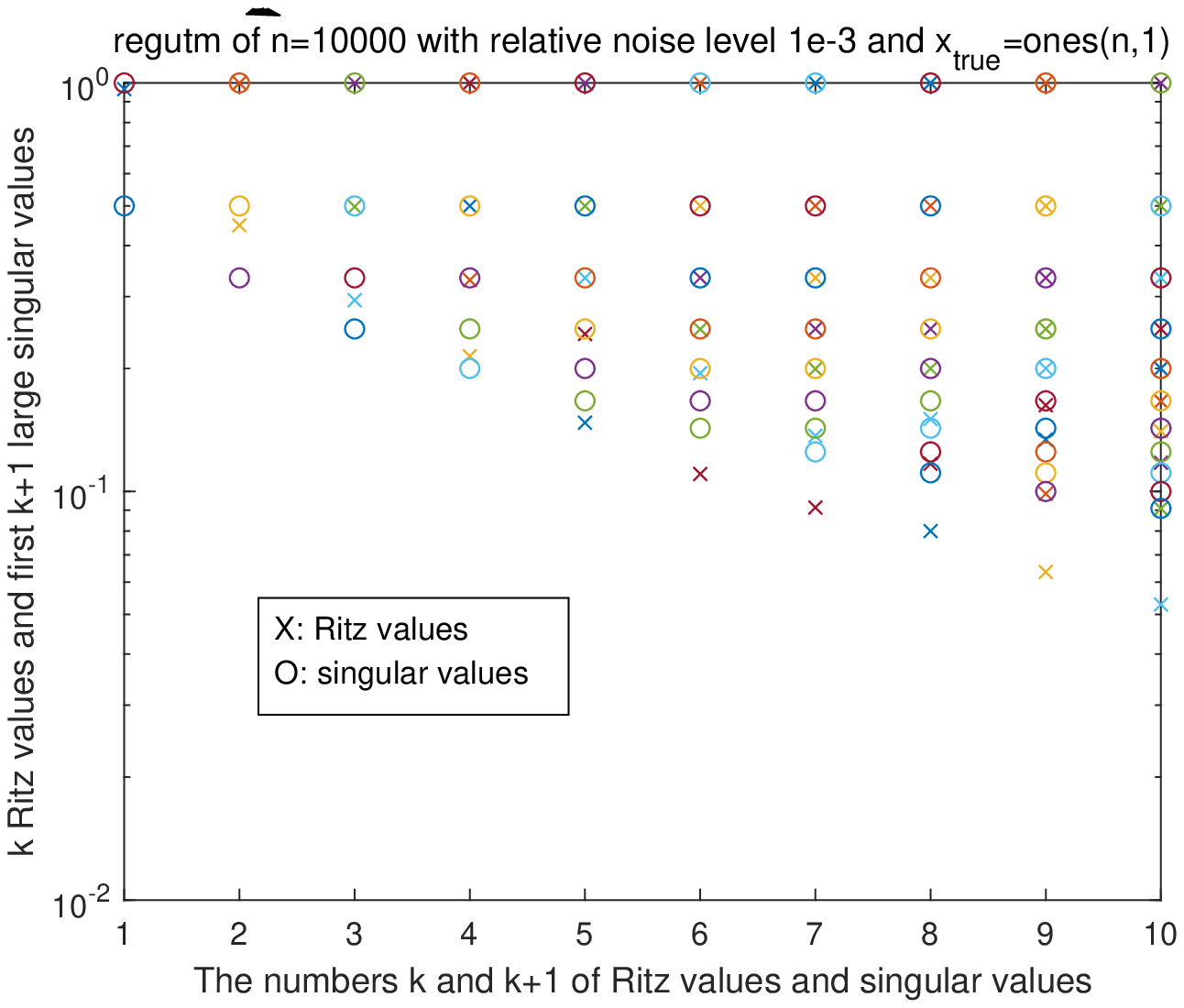}}
  \centerline{(b)}
\end{minipage}
\caption{(a): The accuracy $\gamma_k$ of rank $k$ approximations and
the singular values $\sigma_k=\frac{1}{k}$;
(b) The $k$ Ritz values $\theta_i^{(k)}$ and
the first $k+1$ large singular values $\sigma_i$.}
 \label{fig1}
\end{figure}

\begin{figure}
\begin{minipage}{0.48\linewidth}
  \centerline{\includegraphics[width=6.0cm,height=4.5cm]{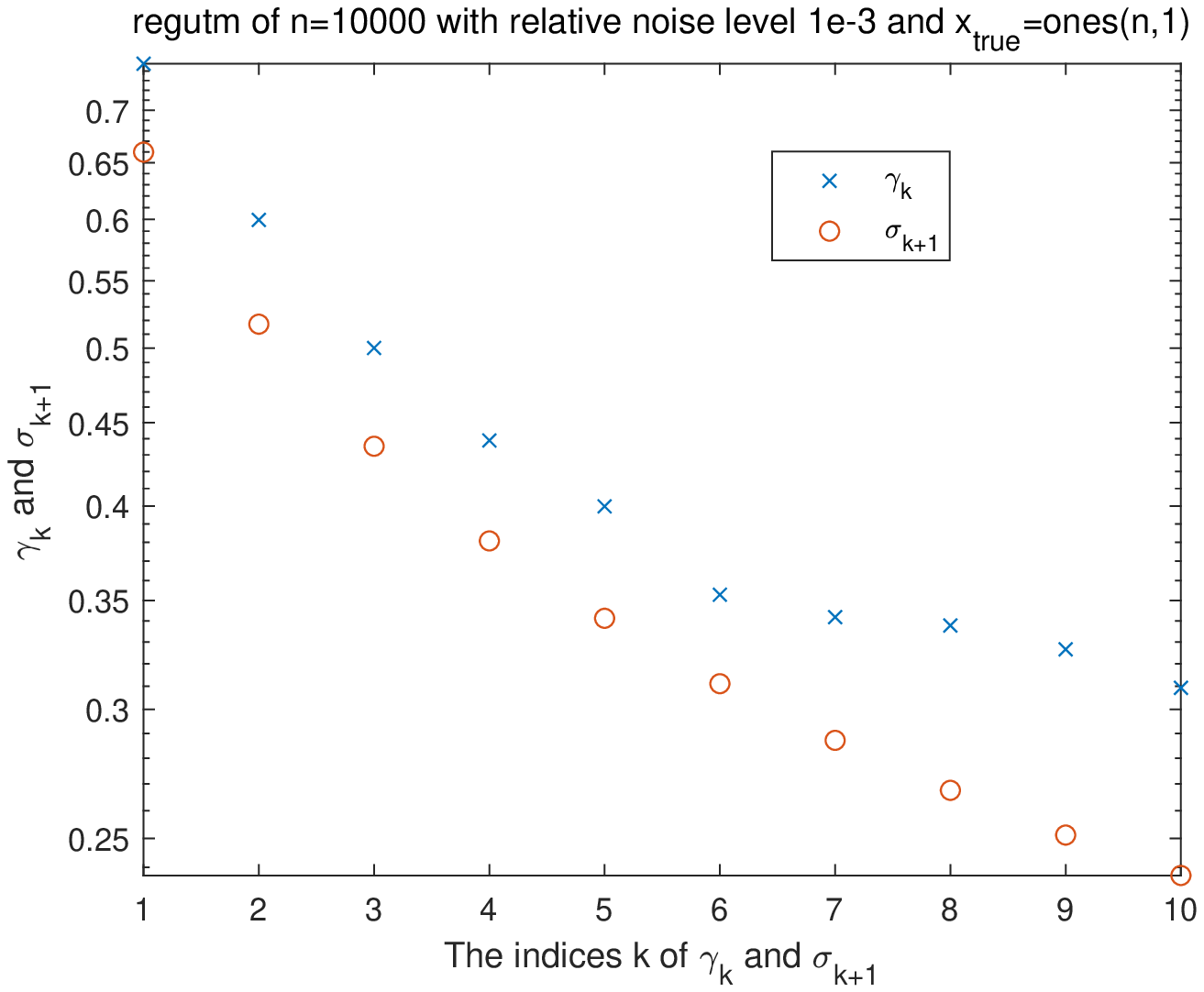}}
  \centerline{(a)}
\end{minipage}
\hfill
\begin{minipage}{0.48\linewidth}
  \centerline{\includegraphics[width=6.0cm,height=4.5cm]{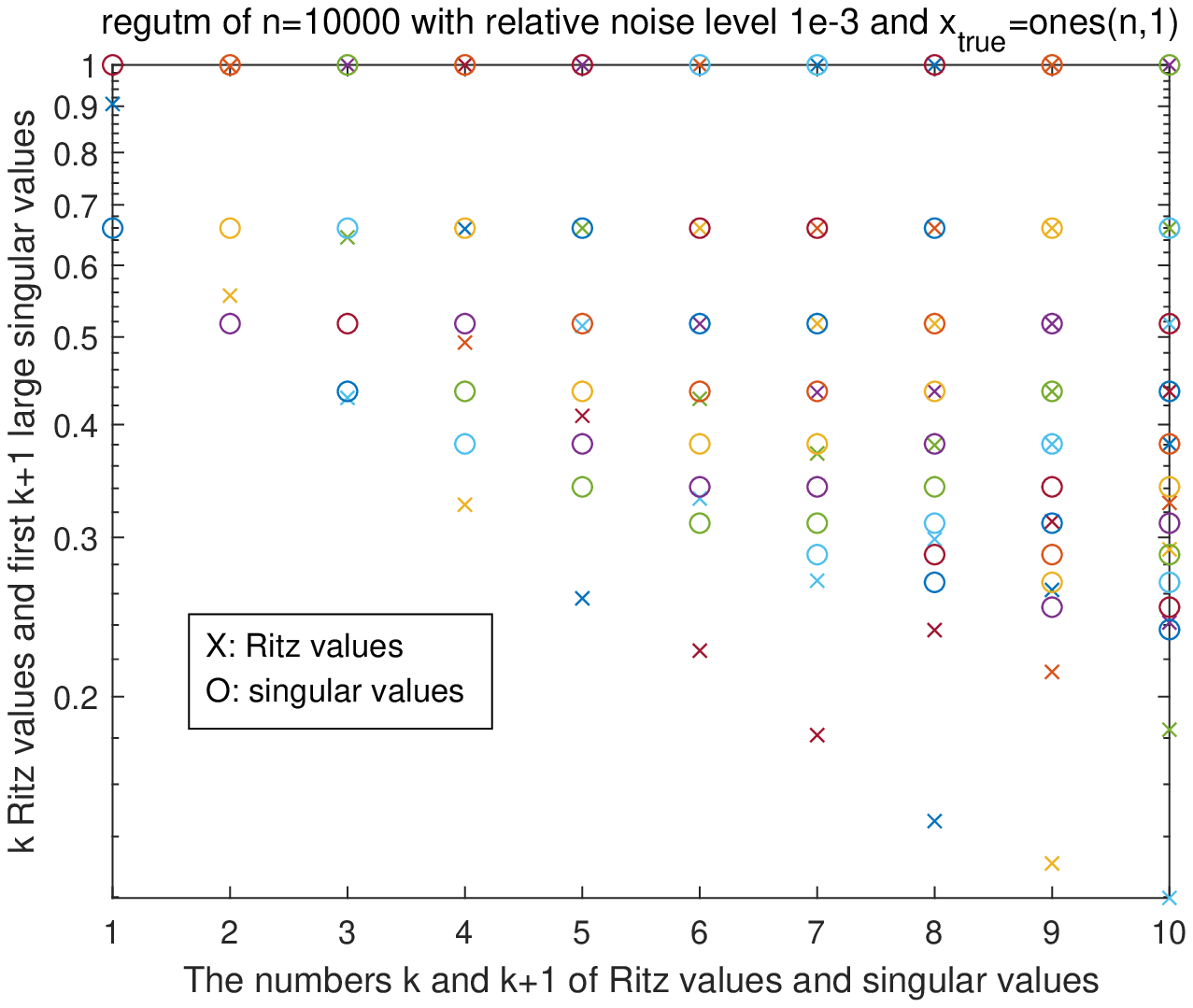}}
  \centerline{(b)}
\end{minipage}
\caption{(a): The accuracy $\gamma_k$ of rank $k$ approximations
and
the singular values $\sigma_k=\frac{1}{k^{0.6}}$;
(b) The $k$ Ritz values $\theta_i^{(k)}$ and
the first $k+1$ large singular values $\sigma_i$.} \label{fig2}
\end{figure}

\begin{figure}
\begin{minipage}{0.48\linewidth}
  \centerline{\includegraphics[width=6.0cm,height=4.5cm]{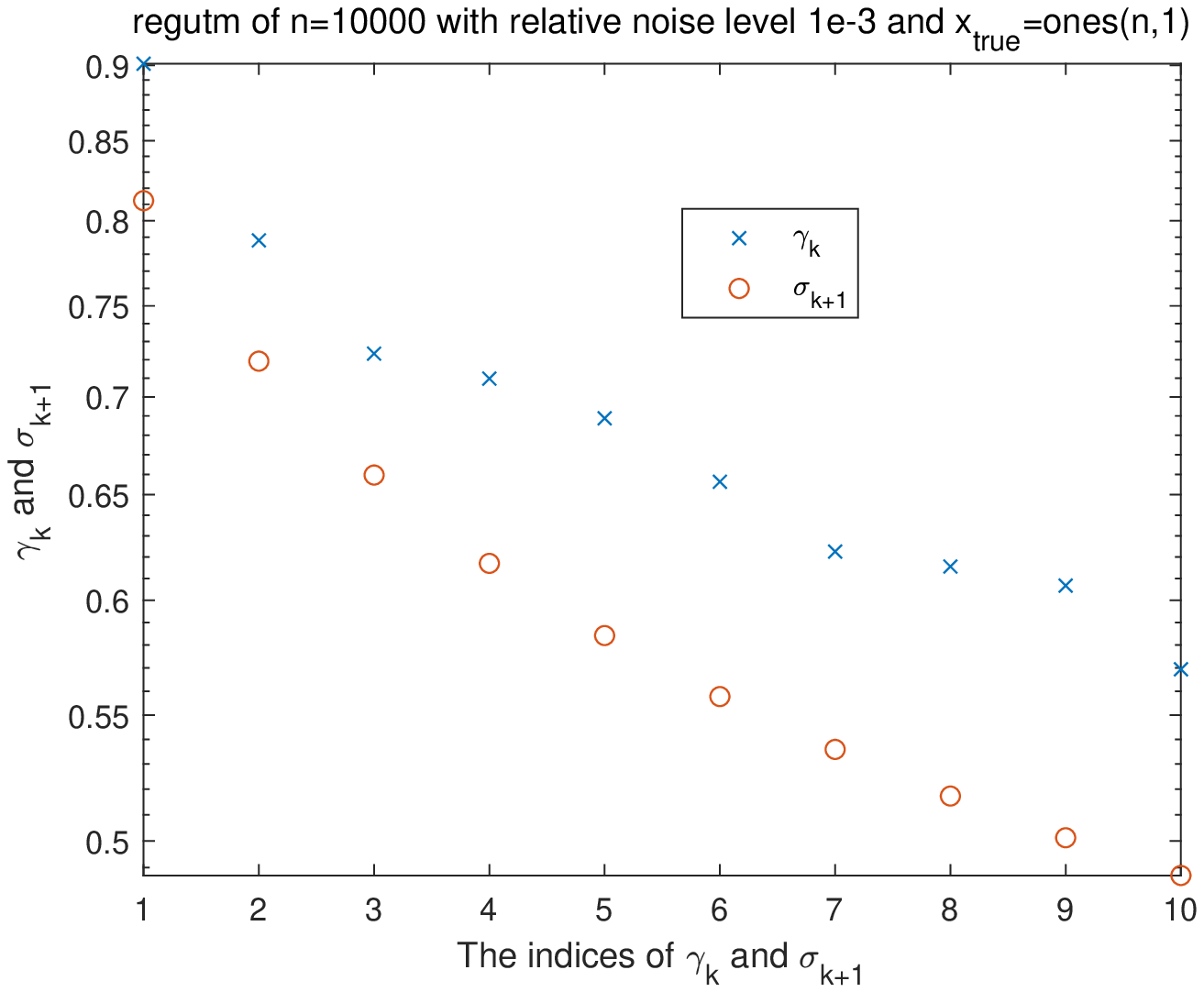}}
  \centerline{(a)}
\end{minipage}
\hfill
\begin{minipage}{0.48\linewidth}
  \centerline{\includegraphics[width=6.0cm,height=4.5cm]{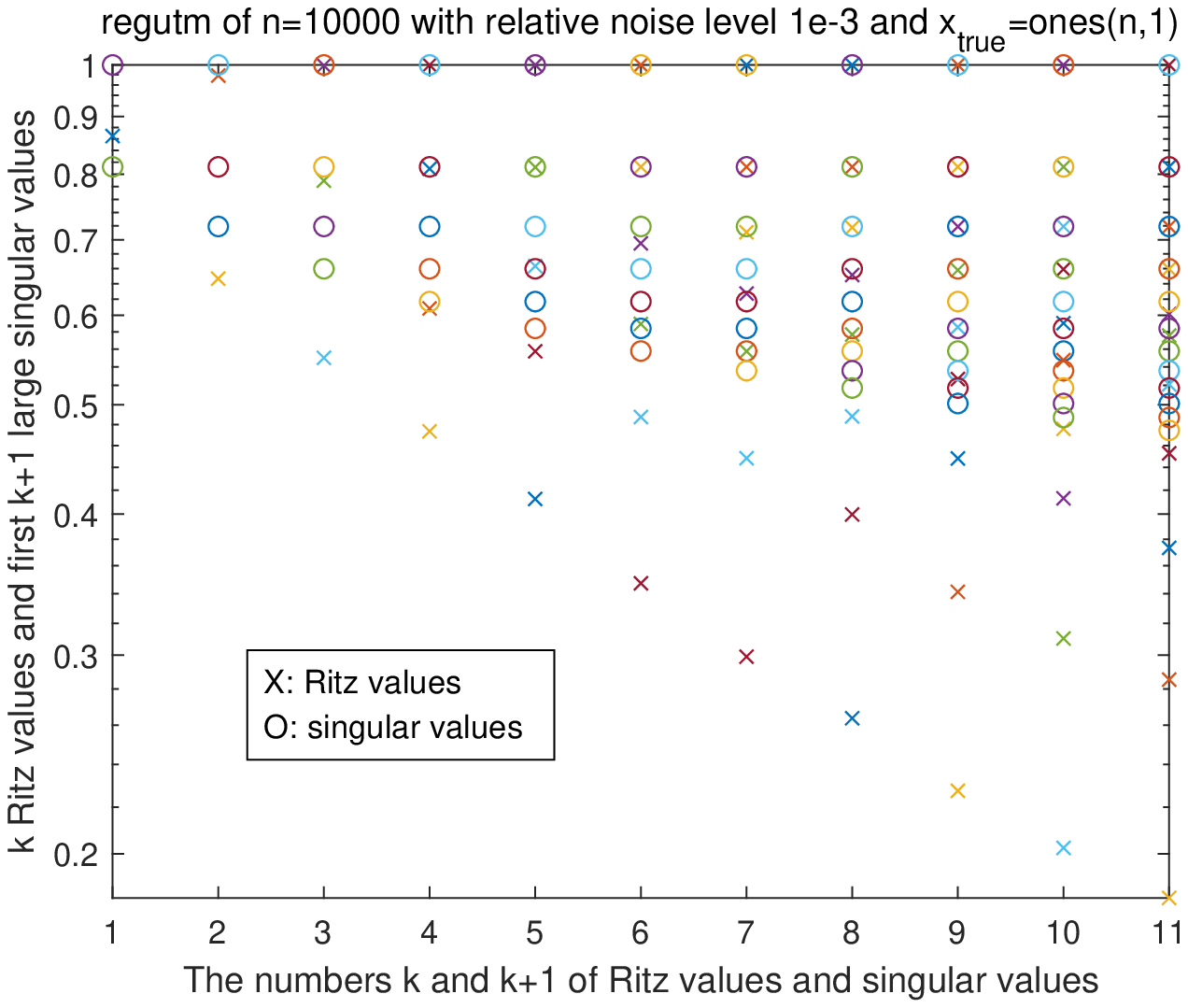}}
  \centerline{(b)}
\end{minipage}
\caption{(a): The accuracy $\gamma_k$ of rank $k$ approximations and
the singular values $\sigma_k=\frac{1}{k^{0.3}}$;
(b) The $k$ Ritz values $\theta_i^{(k)}$ and
the first $k+1$ large singular values $\sigma_i$.}
 \label{fig3}
\end{figure}

Figures~\ref{fig1}--\ref{fig3} (a) draw the accuracy $\gamma_k$
of rank $k$ approximations $P_{k+1}B_kQ_k^T$ and $\sigma_{k+1}$,
and Figures~\ref{fig1}--\ref{fig3} (b) depict the locations of
the $k$ Ritz values $\theta_i^{(k)}$ and the first $k+1$ large singular
values $\sigma_i$ of $A$ for the three $\alpha$'s.
Table~\ref{tab1} highlights these figures
and lists some precise results, where
the first column indicates the location of $\gamma_k$ and to
which singular value it is closer, and the second column ''{\sf near best''}
indicates whether or not $P_{k+1}B_k Q_k^T$ is a
near best rank $k$ approximation to $A$ in the sense of \eqref{near},
the third column ''$\#<\sigma_{k+1}$''
denotes the number of $\theta_i^{(k)}<\sigma_{k+1},\ i=1,2,\ldots,k$,
and the fourth column ''{\sf  $\#$ attained}'' denotes whether
or not the number in the third column attains
the maximum possible number of $\theta_i^{(k)}<\sigma_{k+1},\ i=1,2,\ldots,k$,
which is indicated in Theorem~\ref{disorder}.


\begin{table}[h]
\centering
  \caption{The accuracy $\gamma_k$ of rank $k$ approximations and the
  approximation behavior of the Ritz values $\theta_i^{(k)},\ i=1,2,\ldots,k$.}
  \label{tab1}
  \centering{$\sigma_k=\frac{1}{k},\ k=1,2,\ldots,n$}
  \begin{minipage}[t]{1\textwidth}
  \begin{tabular*}{\linewidth}{lp{2.4cm}p{2.4cm}p{2.4cm}p{2.4cm}p{2.4cm}}
  \toprule[0.6pt]
     The location of $\gamma_k$           &{\sf near best}
     &{\sf $\# <\sigma_{k+1}$}    &{\sf $\#$ attained}  \\ \midrule[0.6pt]
     $\gamma_5\in (\sigma_6,\sigma_5)$ closer
     to $\sigma_6$    &yes     &1      &yes            \\
     $\gamma_6\in (\sigma_7,\sigma_6)$ closer
     to $\sigma_7$    &yes     &1      &yes            \\
$\gamma_7\in (\sigma_7,\sigma_6)$ closer
     to $\sigma_6$    &no     &1      &no            \\
     $\gamma_8\in (\sigma_8,\sigma_7)$ closer
     to $\sigma_7$    &no     &1      &no            \\
     $\gamma_9\in (\sigma_8,\sigma_7)$ closer
     to $\sigma_8$    &no     &2     &no            \\
$\gamma_{10}\in (\sigma_{10},\sigma_9)$ closer
     to $\sigma_{10}$    &no     &2      &yes            \\
    \bottomrule[0.6pt]
     \end{tabular*}\\[2pt]
      \end{minipage}
      \centering{$\sigma_k=\frac{1}{k^{0.6}},\ k=1,2,\ldots,n$}
  \begin{minipage}[t]{1\textwidth}
   \begin{tabular*}{\linewidth}{lp{2.4cm}p{2.4cm}p{2.4cm}p{2.4cm}p{2.4cm}}
  \toprule[0.6pt]
     The location of $\gamma_k$           &{\sf near best}
     &{\sf $\#<\sigma_{k+1}$}    &{\sf  $\#$ attained}  \\ \midrule[0.6pt]
     $\gamma_3\in (\sigma_4,\sigma_3)$ closer
     to $\sigma_3$    &no     &1      &no            \\
     $\gamma_4\in (\sigma_4,\sigma_3)$ closer
     to $\sigma_4$    &no    &1      &no           \\
$\gamma_5\in (\sigma_6,\sigma_5)$ closer
     to $\sigma_5$    &no     &1      &no            \\
$\gamma_6\in (\sigma_{6},\sigma_5)$ closer
     to $\sigma_{6}$    &no     &1      &no            \\
$\gamma_7\in (\sigma_{7},\sigma_6)$ closer
     to $\sigma_{7}$    &no     &2      &yes            \\
     $\gamma_8\in (\sigma_{7},\sigma_6)$ closer
     to $\sigma_{7}$    &no     &2      &no           \\
$\gamma_8\in (\sigma_{7},\sigma_6)$ closer
     to $\sigma_{7}$    &no     &2      &no            \\
$\gamma_9\in (\sigma_{7},\sigma_6)$ closer
     to $\sigma_{7}$    &no     &2      &no          \\
$\gamma_{10}\in (\sigma_{8},\sigma_7)$ closer
     to $\sigma_{7}$    &no     &2      &no           \\
    \bottomrule[0.6pt]
    \end{tabular*}\\[2pt]
    \end{minipage}
  \centering{$\sigma_k=\frac{1}{k^{0.3}},\ k=1,2,\ldots,n$}
  \begin{minipage}[t]{1\textwidth}
 \begin{tabular*}{\linewidth}{lp{2.4cm}p{2.4cm}p{2.4cm}p{2.4cm}p{2.4cm}}
  \toprule[0.6pt]
     The location of $\gamma_k$           &{\sf near best}
     &{\sf $\#<\sigma_{k+1}$}    &{\sf $\#$ attained}  \\ \midrule[0.6pt]
    $\gamma_2\in (\sigma_3,\sigma_2)$ closer
     to $\sigma_2$    &no     &1      &yes            \\
     $\gamma_3\in (\sigma_4,\sigma_3)$ closer
     to $\sigma_3$    &no     &1      &no            \\
     $\gamma_4\in (\sigma_4,\sigma_3)$ closer
     to $\sigma_4$    &no    &2      &yes           \\
$\gamma_5\in (\sigma_6,\sigma_5)$ closer
     to $\sigma_4$    &no     &2      &no            \\
$\gamma_6\in (\sigma_{5},\sigma_4)$ closer
     to $\sigma_{4}$    &no     &2      &no            \\
$\gamma_7\in (\sigma_{5},\sigma_4)$ closer
     to $\sigma_{5}$    &no     &2      &no            \\
     $\gamma_8\in (\sigma_{5},\sigma_4)$ closer
     to $\sigma_{5}$    &no     &3      &no           \\
$\gamma_8\in (\sigma_{5},\sigma_4)$ closer
     to $\sigma_{5}$    &no     &3      &no            \\
$\gamma_9\in (\sigma_{5},\sigma_4)$ closer
     to $\sigma_{7}$    &no     &3      &no          \\
$\gamma_{10}\in (\sigma_{6},\sigma_5)$ closer
     to $\sigma_{6}$    &no     &4      &no           \\
    \bottomrule[0.6pt]
    \end{tabular*}\\[2pt]
  \end{minipage}
\end{table}

Several comments are made in order on the figures and table.

Firstly, for the
test problem with $\alpha=1$, we observe from Figure~\ref{fig1} (a) that
$P_{k+1}B_kQ_k^T$ is a near best rank $k$ approximation to $A$
for $k=1,2,\ldots,6$ and afterwards it is not any longer. However,
starting from $k=5$ onwards, the $\theta_i^{(k)}$ do not approximate
the $k$ large singular values $\sigma_i$ of $A$ in natural order,
as is clearly displayed in Figure~\ref{fig1} (b).
This indicates that the near best rank $k$ approximations
for $k=5,6$ cannot guarantee that the $\theta_i^{(k)}$
approximate the $\sigma_i$ in natural order,
confirming our theory. Furthermore,
it is seen from Table~\ref{tab1} that in each of these two cases there
is exactly one Ritz value $\theta_k^{(k)}<\sigma_{k+1}$, which coincides
with the maximum possible number estimated by Theorem~\ref{disorder}.

Secondly, by comparing (a) with (b) in Figures~\ref{fig1}--\ref{fig3}
correspondingly, we observe that there
is always at least one Ritz value $\theta_k^{(k)}<\sigma_{k+1}$ whenever
$P_{k+1}B_kQ_k^T$ is not a near best rank $k$ approximation to $A$.
We have described these features more clearly in Table~\ref{tab1}.

Thirdly, by inspecting Figures~\ref{fig1}--\ref{fig3} and Table~\ref{tab1}, we
find that, for the same $k$, the smaller $\alpha$ is, the less accurate
the rank $k$ approximation, since for the same
$k$ it is clear that $\gamma_k$ is further away from $\sigma_{k+1}$, which
can be seen from the interval in which $\gamma_k$ lies.
This justifies that it is harder to
generate a good rank $k$ approximation when the decay
of the singular values becomes slower. Particularly,
the rank $k$ approximations are not near
best for $\alpha=0.6$ and $0.3$
from iterations $k=2$ and $k=1$ upwards, respectively.

Fourthly, we observe from Figures~\ref{fig1}--\ref{fig2} (b)
that the smaller $\alpha$ is, the earlier
the Ritz values $\theta_i^{(k)}$ fail to approximate the
$\sigma_i$ in natural order. Precisely, we can see from the figures
that such $k$'s are 5, 3 and 2 for $\alpha=1,0.6 $ and $0.3$, respectively.
Generally, we deduce that the approximations in natural order fail sooner
as the decay of the singular values is slower.

Fifthly, for each $\alpha$, we see from Figures~\ref{fig1}--\ref{fig3}
and Table~\ref{tab1} that as $k$ increases,
the rank $k$ approximation generally becomes poorer
and the number of $\theta_i^{(k)}<\sigma_{k+1}$ exhibits an increasing tendency.

Sixthly, for different $\alpha$, we observe from Table~\ref{tab1}
that for the same $k$ and the smaller $\alpha$, the
number of $\theta_i^{(k)}<\sigma_{k+1},\ i=1,2,\ldots,k$ is
at least nondecreasing and often increases.

Seventhly, for each $\alpha$ and given maximum $k=10$
there is always at least one iteration
$j\leq k$ at which the number of $\theta_i^{(j)}<\sigma_{j+1}$
is exactly equal to its possible maximum; see
the fourth column of Table~\ref{tab1}.

\begin{figure}
\begin{minipage}{0.48\linewidth}
  \centerline{\includegraphics[width=6.0cm,height=4.5cm]{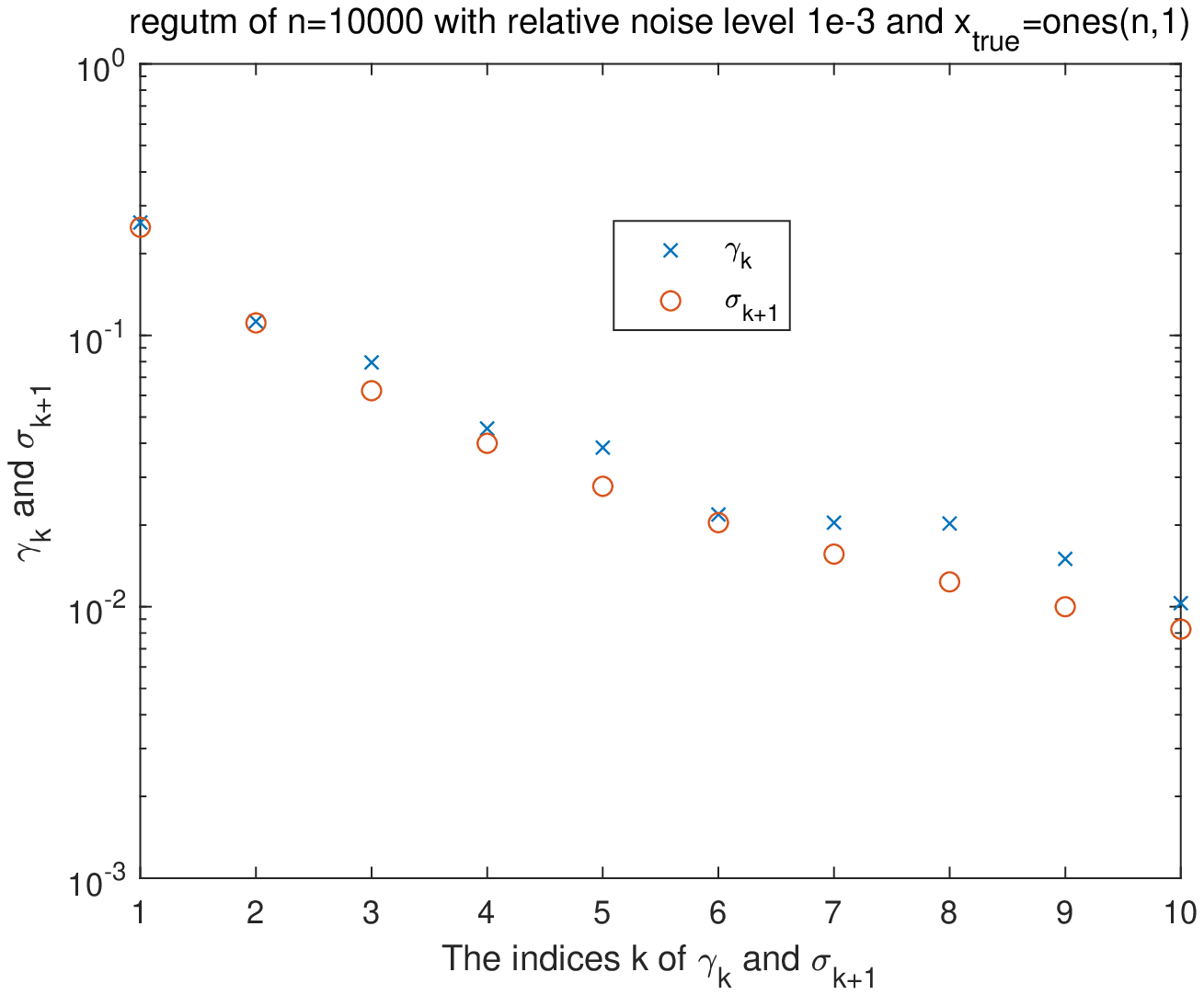}}
  \centerline{(a)}
\end{minipage}
\hfill
\begin{minipage}{0.48\linewidth}
  \centerline{\includegraphics[width=6.0cm,height=4.5cm]{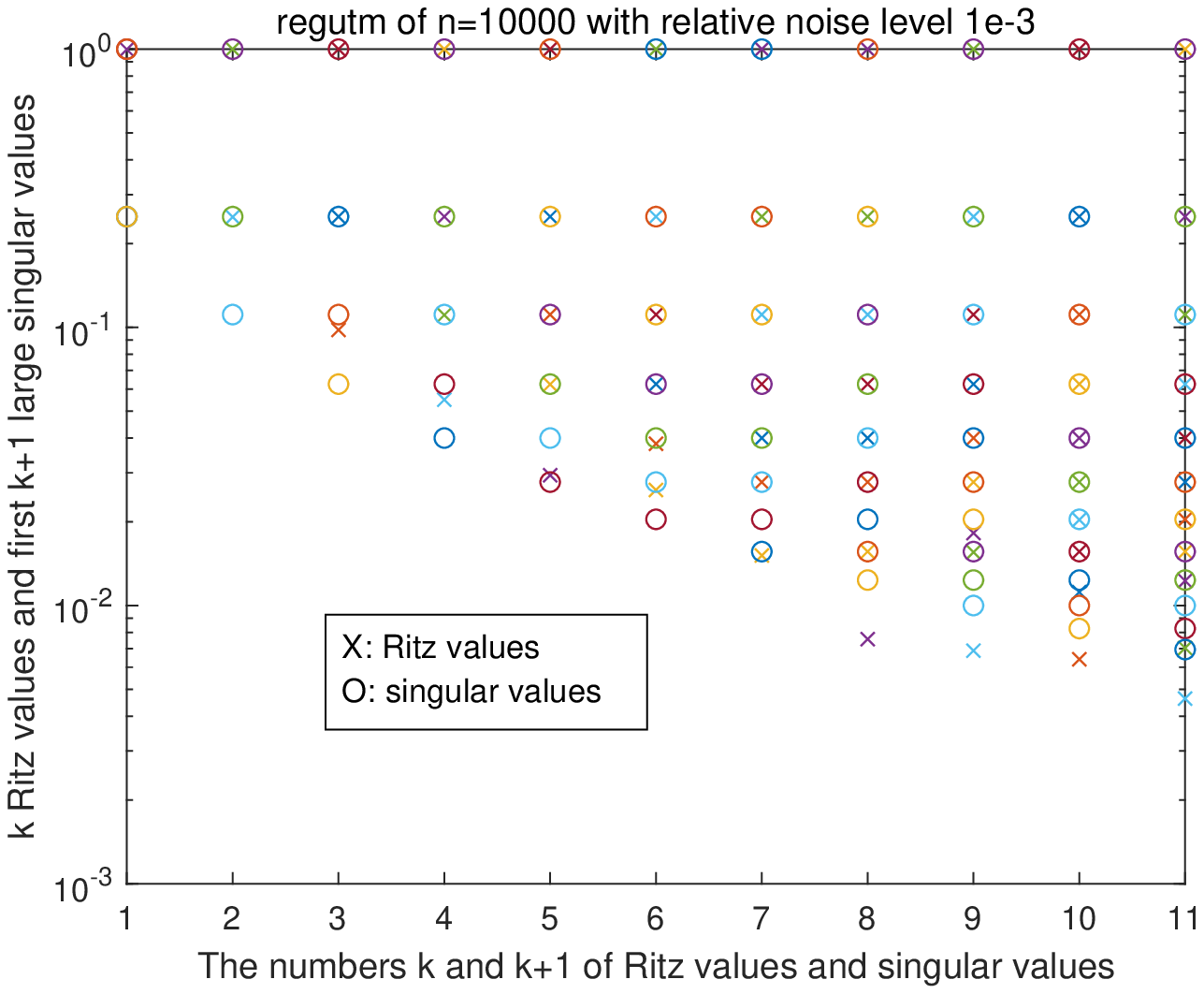}}
  \centerline{(b)}
\end{minipage}
\caption{(a): The accuracy $\gamma_k$ of rank $k$ approximations and
the singular values $\sigma_k=\frac{1}{k^2}$;
(b) The $k$ Ritz values $\theta_i^{(k)}$ and
the first $k+1$ large singular values $\sigma_i$.}
 \label{fig4}
\end{figure}

Next, we report the results on {\sf regutm} of $m=n=10,000$
with $\sigma_k=\frac{1}{k^2}$, i.e., $\alpha=2$,
$x_{true}=ones(n,1)$ and the relative noise level $\varepsilon=10^{-3}$.
This is a moderately ill-posed problem with $\alpha>1$ fairly. We aim to
justify Theorems~\ref{nearapprox}--\ref{ritzvalue} and the second part of
Theorem~\ref{nearappro}. Figure~\ref{fig4} (a)-(b) depict the accuracy
$\gamma_k$ versus $\sigma_{k+1}$ and the $k$ Ritz values $\theta_i^{(k)}$
and the first $k+1$ singular values $\sigma_i,\ k=1,2,\ldots,10$, respectively.

Clearly, we observe from the figure that $P_{k+1}B_kQ_k^T$
is a near best rank $k$ approximation to $A$ until $k=7$ and
afterwards it is not any longer. Correspondingly,
the $k$ Ritz values $\theta_i^{(k)}$ interlace the first $k+1$ large singular
values $\sigma_i$ until the same $k$. Afterwards, such interlacing
property is lost, and the Ritz values do not approximate the singular
values in natural order any more. These justify
Theorems~\ref{nearapprox}--\ref{ritzvalue}, in which, for a
fixed $\alpha>1$, the
sufficient conditions \eqref{condition1} and \eqref{condm}
fail to meet when $k$ increases up to some point, since the
left-hand sides strictly increase and the right-hand sides strictly
decrease with $k$ for the fixed $\alpha$. They also
confirm the second part of Theorem~\ref{nearappro}, which requires
that $\alpha>1$ sufficiently for a given not small $k$,
while $\alpha=2$ cannot meet this requirement for $k>7$.

By comparing Figure~\ref{fig4} with Figures~\ref{fig1}--\ref{fig3},
we mention that $k=7$ is
considerably bigger than those for $\alpha=1,0.6 $ and
$0.3$, after which $P_{k+1}B_kQ_k^T$ is not a near best rank
$k$ approximation and the Ritz values fail to approximate the
singular values in natural order.
This again confirms that, for the same $k$, the rank $k$ approximation
$P_{k+1}B_kQ_k^T$ is more accurate and the Ritz values are more likely
to approximate the singular values in natural order for a bigger $\alpha$.

\begin{figure}
\begin{minipage}{0.48\linewidth}
  \centerline{\includegraphics[width=6.0cm,height=4.5cm]{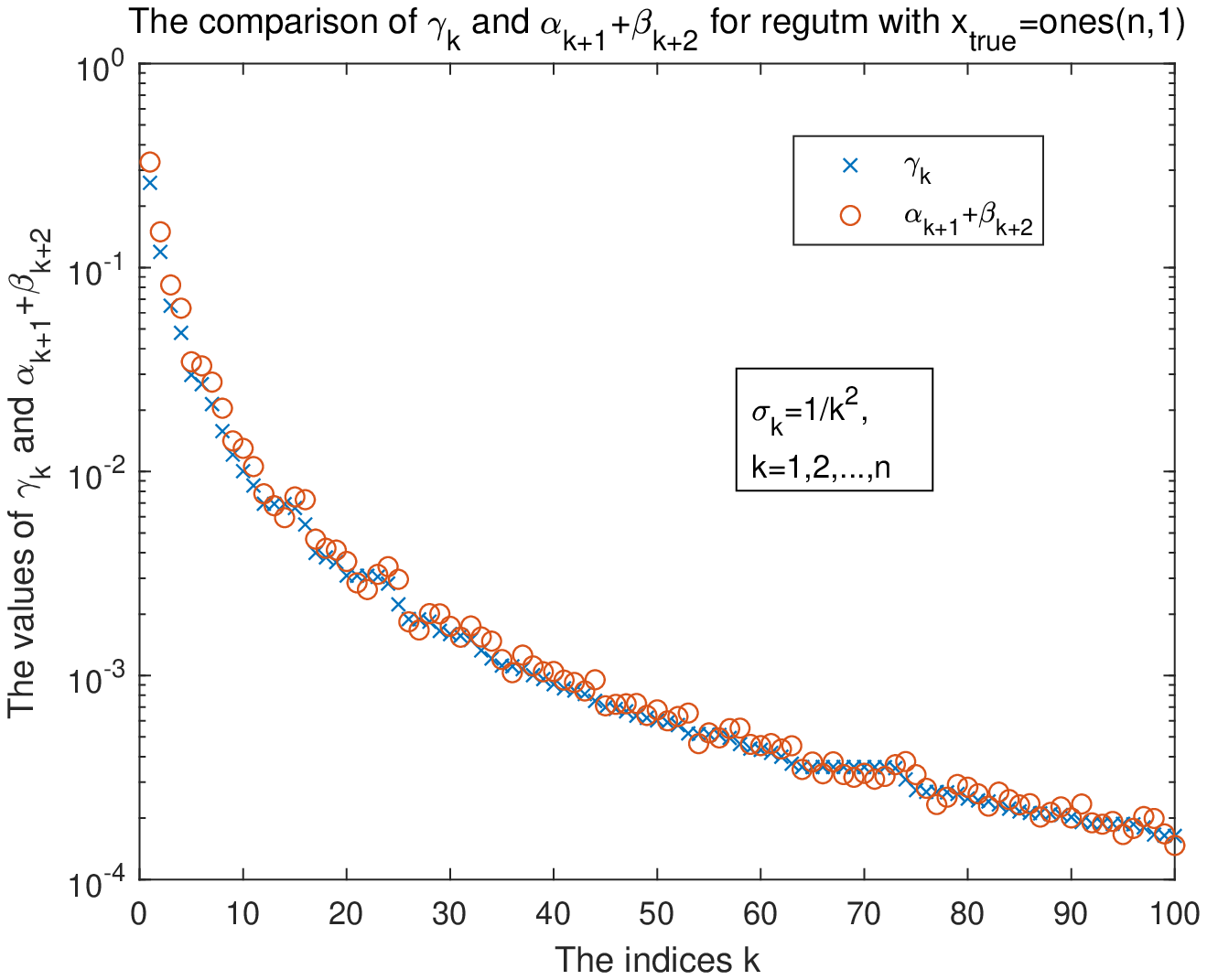}}
  \centerline{(a)}
\end{minipage}
\hfill
\begin{minipage}{0.48\linewidth}
  \centerline{\includegraphics[width=6.0cm,height=4.5cm]{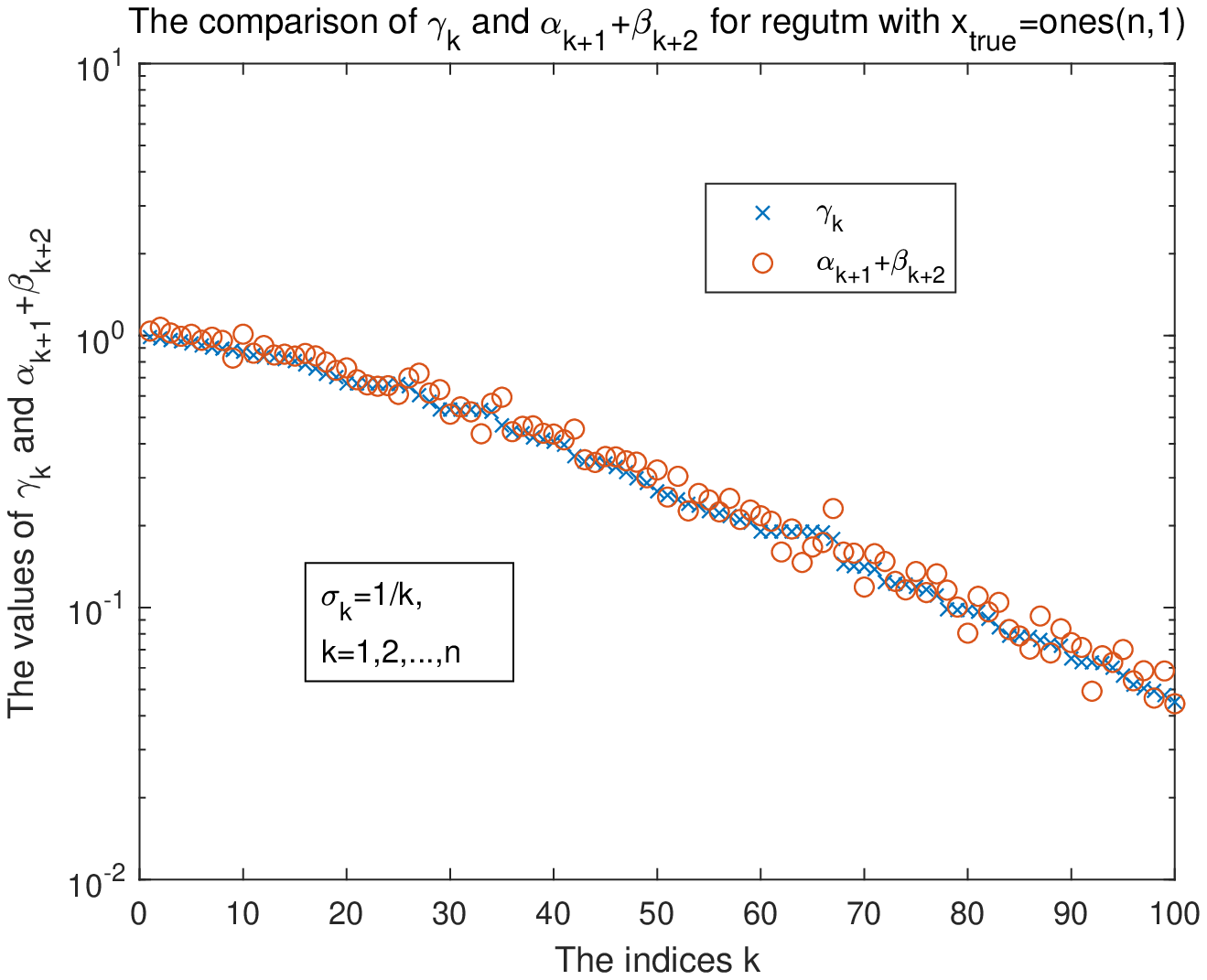}}
  \centerline{(b)}
\end{minipage}
\begin{minipage}{0.48\linewidth}
  \centerline{\includegraphics[width=6.0cm,height=4.5cm]{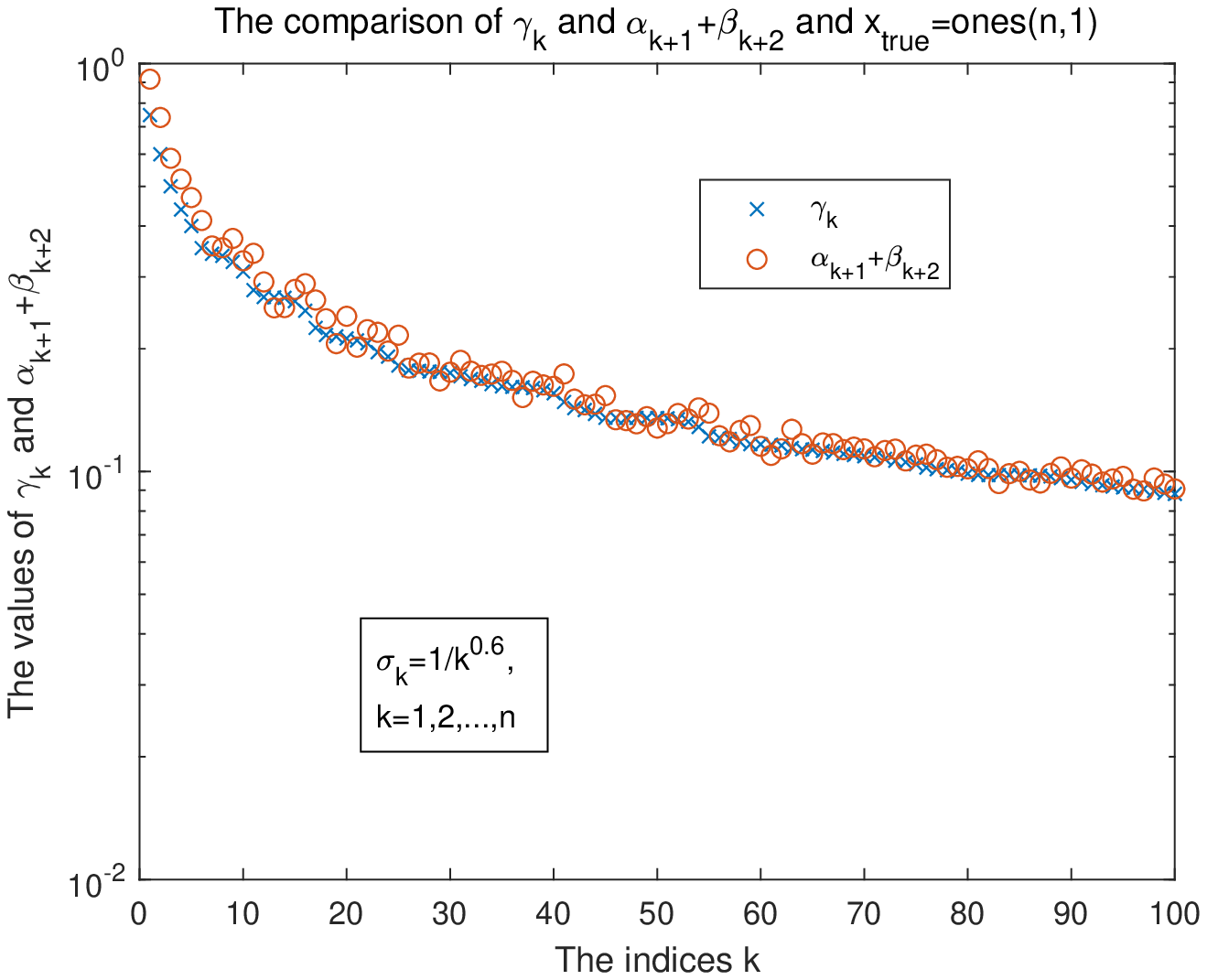}}
  \centerline{(c)}
\end{minipage}
\hfill
\begin{minipage}{0.48\linewidth}
  \centerline{\includegraphics[width=6.0cm,height=4.5cm]{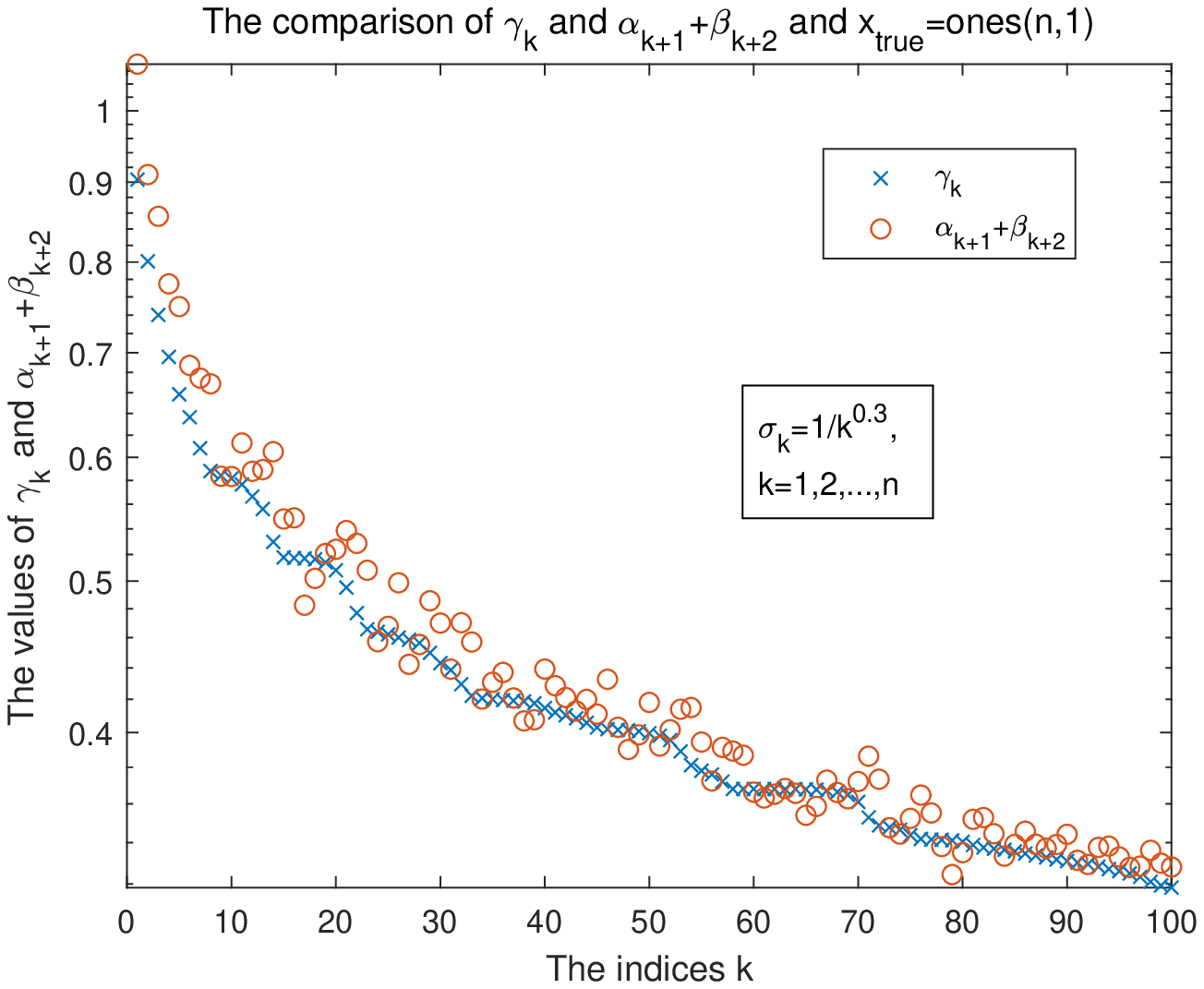}}
  \centerline{(d)}
\end{minipage}
\caption{$\gamma_k$ versus $\alpha_{k+1}+\beta_{k+2}$ for {\sf regutm} of $m=n=10,000$
with $\alpha=2,1,0.6$ and $0.3$.}
 \label{fig5}
\end{figure}

Finally, we justify Theorem~\ref{main2} and the first two remarks followed.
Figure~\ref{fig5} depicts $\gamma_k$ versus $\alpha_{k+1}+\beta_{k+2}$
for {\sf regutm} with $\alpha=2,1,0.6$ and $0.3$. Clearly, we see from the
figure that $\alpha_{k+1}+\beta_{k+2}$
decays exactly as fast as $\gamma_k$. Therefore, independent of the
degree of ill-posedness, we can reliably judge the decreasing property
and tendency of the practically uncomputable accuracy $\gamma_k$
by the available $\alpha_{k+1}+\beta_{k+2}$ with little cost.
In addition, it is clearly seen from the figure that $\gamma_k$ decays
faster for a bigger $\alpha$ than for a smaller $\alpha$.

\section{Conclusions}\label{concl}

For the large-scale ill-posed problem \eqref{eq1}, LSQR
is a most commonly used Krylov solver for general purposes.
It has general regularizing effects and exhibits semi-convergence.
If LSQR has already found
best possible 2-norm filtering regularized solutions, then
it has the full regularization. In this case,
we simply stop it after a few iterations of semi-convergence,
and complicated hybrid variants are not needed. The semi-convergence
of LSQR, in principle,
can be determined by a suitable parameter-choice method, such as
the L-curve criterion and the discrepancy principle \cite{hansen98,hansen10}

In the simple singular value case, the author in \cite{jia18a,jia18b}
has proved that, for the severely and moderately ill-posed problems
with suitable $\rho>1$ and $\alpha>1$,
the Ritz values $\theta_i^{(k)}$ approximate
the large singular values of $A$ in natural order until the semi-convergence,
so that LSQR has the full regularization.
On the other hand, however, if $\alpha>1$ not enough, the approximation in this
order cannot be ensured; if $\alpha\leq 1$, such approximation property
may hold only for $k$ very small. In this paper, we have nontrivially
extended the results in \cite{jia18a,jia18b} to the multiple singular value
case, and drawn the same conclusions.

As a major contribution, we have made an
in-depth analysis on best, near best and general rank $k$ approximations
to $A$ for more general ill-posed problems with $\alpha>0$, which include
mildly ill-posed problems, and derived some insightful properties
of the Ritz values $\theta_i^{(k)},\ i=1,2,\ldots,k$.
Our results have shown that
general best or near best rank $k$ approximations do not guarantee that
$\theta_i^{(k)},\ i=1,2,\ldots,k$ approximate
large singular values of $A$ in natural order for $0<\alpha\leq 1$.
We have proved that, for the same $k$, the smaller $\alpha$ is,
the less accurate the rank $k$ approximation is, and the more
Ritz values smaller than $\sigma_{k+1}$. Numerical experiments
have confirmed the theoretical results.
These results apply to the other Krylov solvers CGME, MINRES, MR-II, GMRES
and RRGMRES as well.



\bibliographystyle{amsplain}


\end{document}